\documentclass[12pt]{article}

\usepackage{latexsym}
\usepackage{graphicx}
\usepackage[export]{adjustbox}
\usepackage{amssymb,amsmath}
\usepackage[english]{babel}
\usepackage{color}

\newcommand{\Gn}{\ensuremath{\boldsymbol \Gamma_{\!\rm{n}}}} 
\newcommand{\Gpr}{\ensuremath{\boldsymbol \Gamma_{\!\rm{pr}}}}
\newcommand{\Gpo}{\ensuremath{\boldsymbol \Gamma_{\!\rm{pt}}}}
\newcommand{\dn}{\ensuremath{\mathbf  d}} 
\newcommand{\dodd}{\ensuremath{\mathbf d_{\rm odd}}}
\newcommand{\deven}{\ensuremath{\mathbf d_{\rm even}}}
\newcommand{\nnu}{\ensuremath{\boldsymbol \nu}}
\newcommand{\nnumap}{\ensuremath{\nnu_{\!\rm{\scriptscriptstyle MAP}}}}

\begin{document}

\title{Object based Bayesian full-waveform inversion for shear elastography}

\author{Ana Carpio (Universidad Complutense de Madrid),\\
Elena Cebri\'an (Universidad de Burgos),  \\
Andrea Guti\'errez (Universidad Complutense de Madrid)}

\maketitle

{\bf Abstract.} 
We develop a computational framework to quantify uncertainty in
shear elastography imaging of anomalies in tissues. We adopt a Bayesian
inference formulation. Given the observed data, a forward model and their 
uncertainties, we find the posterior probability of parameter fields
representing the geometry of the anomalies and their shear moduli. 
To construct a prior probability, we exploit the topological energies of 
associated objective functions. We demonstrate the approach on synthetic 
two dimensional tests with smooth and irregular shapes.
Sampling the posterior distribution by Markov Chain 
Monte Carlo (MCMC) techniques we obtain statistical information on the 
shear moduli and the geometrical properties of the anomalies. General 
affine-invariant ensemble  MCMC samplers are adequate for shapes 
characterized by parameter sets of low to moderate dimension. 
However, MCMC methods are computationally expensive. For simple 
shapes, we devise a fast optimization scheme to calculate the maximum a 
posteriori (MAP) estimate representing the most likely parameter values. 
Then, we approximate the posterior distribution by a Gaussian distribution 
found by linearization about the MAP point to capture the main mode at 
a low computational cost.

\section{Introduction}
\label{sec:intro}

Medical imaging is a part of biological imaging that aims to reveal internal 
structures hidden in tissues by non invasive techniques \cite{medical}, 
such as X-ray radiography, magnetic resonance imaging, tomography, 
echography, ultrasound, endoscopy, elastography, tactile imaging, thermography, 
nuclear medicine and holography. From the mathematical point of view, 
they all pose inverse problems that aim to deduce the properties of living 
tissues from observed signals. The typical framework is as follows. A set 
of emitters launch waves which interact with the tissue. The resulting wave 
field  is recorded at a grid of receptors and analyzed to infer the structure 
of the medium. Different imaging techniques differ  in the waves employed 
(electromagnetic, acoustic, thermal, elastic, etc), the arrangement of emitters 
and receivers, and the medium properties monitored. 

Harmless modalities using light \cite{light}, sound \cite{sound} or elastic 
\cite{elastic} beams are particularly interesting due to the absence of 
secondary effects. Elastography is a relatively new imaging technique 
that maps the elastic properties of soft tissue \cite{elastic, tumorstiffness}. 
Cancerous tumors will often be stiffer than the surrounding tissue (prostate 
and breast tumors, for instance), whereas damaged livers are harder than 
healthy ones \cite{tumorspeed, tumorelasticity, tumorstiffness}. 
While existing technology \cite{transient} can distinguish healthy from 
unhealthy  tissue in specific  situations, the study of tissues containing 
multiple anomalies, tiny tumors or little contrast regions may benefit 
from the development of more refined mathematical approaches.  

Here we develop an object based Bayesian full-waveform inversion framework 
for soft tissue shear elastography with topological priors. Instead of tracking 
spatial variations of the elastic constants or the wave speeds within the 
tissue (as is often done in many geophysical and medical applications, 
see \cite{tumorspeed, fichtner_adjoint,  fichtner_bfwi, georg_mcmc, elastic,  
georg_marmousi}, for instance and references therein), we take an inverse 
scattering approach \cite{colton}  and   represent localized anomalies in 
the tissue, such as tumors or fibromas, by  objects with distinctive elastic 
constants  immersed in the background tissue \cite{guzina1}. 
A first advantage of this approach is that anomalies are characterized by 
a few unknowns defining their parametrization and their elastic properties, 
which reduces the computational cost. Moreover, studies in other imaging 
set-ups \cite{keller} suggest that localized inhomogeneities may be more 
precisely captured by looking for abrupt interfaces defining their boundaries, 
and for material parameter variations within them, than by tracking the 
spatial variations of material parameter fields everywhere.
A second advantage is that we can define misfit functionals in terms of object
 shapes and then use the associated topological energies to construct 
 sharp priors at low cost. Because of their ability to suppress oscillations  
 in configurations  with multiple objects, topological energies are used in 
 deterministic  inverse scattering frameworks to find first guesses of  scatterers 
 in nondestructuve materials  testing  \cite{tenergy0, tenergy1} and in 
 biological  applications \cite{tenergybio}.

Depending on the expected complexity, one can choose different representations 
for the boundaries of the anomalies \cite{hohage3d, level_set, fichtner_object, 
palafox}.
Here, we consider two types of star-shaped parametrizations that differ in the 
way the radius is parameterized. The boundary of star-shaped objects is defined 
by an angle dependent distance function (the radius) along rays in all space 
directions \cite{ghattas_radius}. We can reproduce smooth shapes approximating 
the radius by trigonometric polynomials involving just a few parameters 
\cite{hohage2d}. 
Rougher boundaries are better described by high dimensional radius functions 
 \cite{matt_thesis}. 
Both situations are of interest to study anomalies in tissues. Tumors, for 
instance, can display smooth or irregular contours depending on their stage 
and nature.
We  show that we can infer the structure of anomalies in tissues with 
quantified uncertainty by Markov Chain Monte Carlo (MCMC) sampling of 
posterior distributions which use priors constructed by topological energy 
methods and likelihoods defined in terms of the difference of the recorded 
elastography data and synthetic observations generated numerically for 
arbitrary anomalies. 
Affine invariant ensemble samplers \cite{goodmanweare,matt_sampler} work 
reasonably well when we can approximate the radius function by combinations 
of trigonometric polynomials. We can also extract basic information on irregular 
objects defined by higher dimensional radius functions. However, MCMC 
methods are computationally expensive. For simple shapes, we have also 
developed fast methods which first optimize to calculate a maximum a 
posteriori \cite{georg_linearized, sergei} approximation to the anomaly 
parameters  and then sample a linearized approximation of the posterior 
distribution. The computational cost is much lower, but details on the structure
of the posterior distributions, such as multimodality, may be lost.

The paper is organized as follows.
Section \ref{sec:set-up} describes the mathematical model for shear wave 
imaging in the tissue. Section \ref{sec:inverse} formulates the Bayesian 
inversion framework.
Section \ref{sec:pgeometry} explains how to construct priors for the number 
of anomalies and their shapes.  Section \ref{sec:MCMC}  uses ensemble 
MCMC samplers to solve the Bayesian inverse problem and quantify 
uncertainty in the solution for relevant configurations characterized by low 
dimensional parameter sets. Well defined maximum a posteriori (MAP)
approximations are identified. Section \ref{sec:linearized} presents a low  
cost approach which combines optimization to calculate the MAP point and 
linearization of the posterior probability about it to quantify uncertainty.
Finally, section \ref{sec:irregular} adapts affine-invariant ensamble sampling 
methods to infer the structure of high dimensional irregular shapes. A final 
Appendix contains  details on the numerical schemes employed and
parameter choices made.
Section \ref{sec:conclusions} presents our conclusions.

\begin{figure}[!hbt]
\centering
\includegraphics[width=7cm]{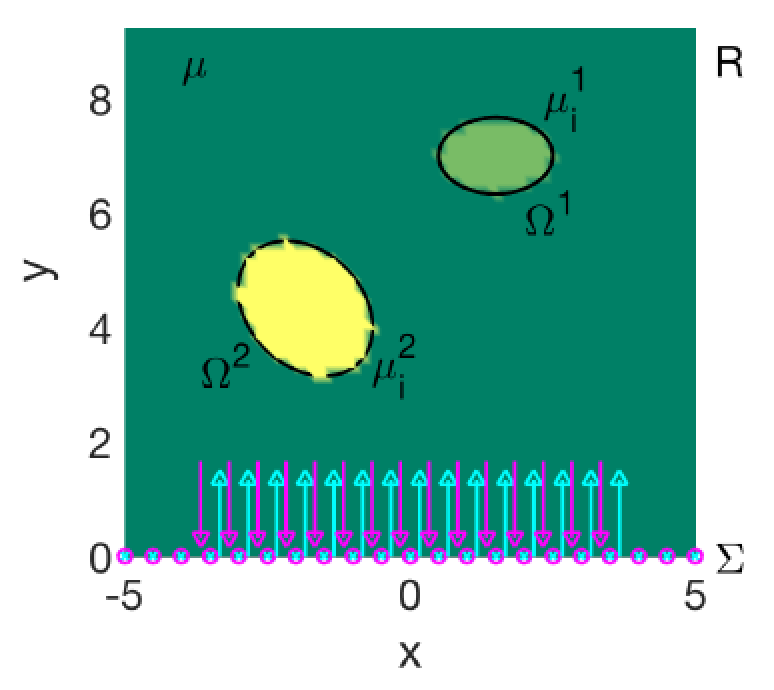} 
\caption{Schematic dimensionless representation of the imaging set-up. 
The emitters (cyan) generate waves which interact with the medium. The 
reflected waves are recorded at the receivers (magenta). Emitters and
receivers are transducers located at the same position. }
\label{fig1}
\end{figure}

\section{Physical set-up}
\label{sec:set-up}

Shear elastography tracks variations in the shear modulus $\mu$, 
which is the property varying  more abruptly from healthy to unhealthy 
tissue, by means of shear waves. 
Elastic waves in a medium split in shear components (shear S-waves) 
and compression components (longitudinal P-waves) \cite{landau}. 
Shear waves are adequate for the depths considered in tissues 
since P-waves travel faster and reach deeper very fast. Moreover, at 
low frequencies, shear waves  are not really  affected by attenuation 
effects in tissues  \cite{transient} and are governed by standard
wave equations.

The imaging set-up is represented in Figure \ref{fig1}. We consider
a medium $R\subset \mathbb R^2$ (the tissue) containing a set of 
anomalies $\Omega = \cup_{\ell=1}^L \Omega^\ell$
and locate a set of emitters  $\mathbf x_j$, $j=1,\ldots,J,$ on a part
$\Sigma$  of the boundary $\partial R$. The medium has density 
$\rho$ and elastic constants $\mu$ and $\lambda$, while the anomalies 
have  density $\rho_{\rm i}$ and elastic constants $\mu_{\rm i}$ and
$\lambda_{\rm i}$. 
The emitted waves  interact with the medium and the resulting wave 
field is recorded at a grid of receivers $\mathbf r_k$, $k=1, \ldots, K$.
Emitter and receivers occupy the same region. They can be interspaced,
or, in some set-ups, overlap. Here, we will consider they are transducers
located at same positions, playing both roles alternatively. Let us
formulate the forward problem that governs the dynamics of the
wave field in this framework.

To simplify, we consider that the waves emitted by the sources are 
governed by the scalar wave equation
\begin{eqnarray}
\begin{array}{ll}  
\rho u_{tt}  - {\rm div } (\mu \nabla u) = f(t) g(\mathbf x),  
& \mathbf x \in R, \; t>0, \\[1ex]
u(\mathbf x, 0) = 0, u_t(\mathbf x, 0) = 0, & \mathbf x \in R,
\end{array} \label{forward_dim}
\end{eqnarray}
where
\begin{eqnarray*}  
\rho(\mathbf x) = \left\{ \begin{array}{ll} 
\rho,    & \mathbf x \in R \setminus \overline{\Omega}, \\
\rho_{\rm i}^\ell, & \mathbf x \in  \Omega^\ell, \; \ell = 1, \ldots, L,
\end{array}  \right.  \label{coefrho}  \\
\mu(\mathbf x) = \left\{ \begin{array}{ll} 
\mu,    & \mathbf x \in R \setminus \overline{\Omega}, \\
\mu_{\rm i}^\ell, & \mathbf x \in  \Omega^\ell, \ell = 1, \ldots, L,
\end{array} \right.
\label{coefmu}
\end{eqnarray*}
with local wave speed   $\sqrt{\mu \over \rho}$ in the healthy tissue 
and $ \sqrt{\mu_{\rm i}^\ell \over \rho_{\rm i}^\ell}$ inside 
each anomaly $\Omega^\ell.$ In tissues,  we have $\rho_{\rm i} \sim \rho$.
We assume that the emitters $\mathbf x_j$, $j=1,\ldots,J,$ 
induce source terms of the form  $f(t) g_j(\mathbf x - \mathbf x_j)$, 
where $g_j$ are smooth functions of narrow support about $0$ that we
sum to obtain $g(\mathbf x)$.  We represent the function $f(t)$ by a Ricker 
wavelet $f(t)= f_0 (1- 2 \pi^2 f_M^2 t^2) e^{-\pi^2 f_M^2 t^2}$ with peak 
frequency $f_M$. The time it takes to move from the initial positive 
maximum to the negative minimum is $T_D = {\sqrt{6} \over 2 \pi f_M}.$ 
After that it approaches zero. 

A whole organ $R$ can be represented by a domain with zero normal 
derivative at  its physical boundary $\partial R$
\begin{eqnarray}  
\displaystyle {\partial u \over \partial \mathbf n} = 0 \quad  \mbox{\rm on } 
\,\,  \partial R.
\label{bcneuman}
\end{eqnarray}
Equations (\ref{forward_dim})-(\ref{bcneuman}) define the forward 
model,  where $\rho \in L^\infty(R)$, $\mu \in L^\infty(R)$,
$\rho \geq \rho_0 >0$ and $\mu \geq \mu_0 >0$. 
Assuming that $R$ and $\Omega$ have $C^1$ boundaries, problem
(\ref{forward_dim})-(\ref{bcneuman}) has a unique solution 
$u \in C([0,\tau];H^1(R))$, $u_t \in C([0,\tau];L^2(R))$, 
$u_{tt} \in L^2(0,\tau;(H^{1}(R))')$, for any  $\tau>0$, see \cite{lions}.  
Here, $H^1(R)$ represents the standard Sobolev space and 
$(H^{1}(R))'$ its dual space \cite{brezis}. 
Since \ $u_{tt}(\mathbf x, 0)=  \in L^2(R)$,
we also have $u_t \in C([0,\tau];H^1(R))$ and 
$u \in C([0,\tau];H^2(R\setminus \overline \Omega))$, see Appendix.
Then, $u(t)$ is defined on  $\Sigma$ and at the receiving sites both in the 
sense of $L^2(\Sigma)$ traces and pointwise \cite{adams,necas}.

Assume that $\Omega_{\rm true}$ represents the true anomalies and 
$\rho_{\rm i,\rm true}$, $\mu_{\rm i,\rm true}$ represent their true 
material  properties. In principle, the values recorded at the receivers 
constitute the data, that is,
$d_{k,\rm true}^m = u(r_k,0,t_m) $, where $u$ is the solution of 
(\ref{forward_dim})-(\ref{bcneuman}) when $\Omega = \Omega_{\rm true}.$
In practice, the recorded data $d_k^m$ are corrupted by different sources 
of noise, that is, $d_k^m = d_{k,\rm true}^m + {\rm noise}.$ We will assume 
that  the additive noise is distributed as a multivariate Gaussian 
${\cal N}(0,\Gn)$ with mean zero and covariance matrix $\Gn$:
\begin{eqnarray}
d_k^m = d_{k,\rm true}^m +  \varepsilon_k^m, \label{noisy}
\end{eqnarray}
for $k=1,\ldots, K$, $m=1,\ldots, M$, where $\varepsilon$ is distributed
according to ${\cal N}(0,\Gn)$.
We consider the noise level for each receiver to be equal and
uncorrelated, so that $\Gn$ is the identity matrix of dimension
$N=KM$ multiplied by $\sigma_{\rm noise}^2.$

\section{Inverse problem}
\label{sec:inverse}

The inverse problem consists in finding the anomalies $\Omega$ and their  
material coefficients $\rho_{\rm i}$ and $\mu_{\rm i}$ such that the solution 
of  the forward problem  agrees, in a way to be specified, 
with  the recorded data. For shear elastography in tissues, we take 
$\rho_{\rm i} \sim \rho$, thus we only have to identify $\mu_{\rm i}$.
In a deterministic framework, one typically resorts to optimization formulations: 
Find objects $\Omega$ and parameters $\mu_{\rm i}$ minimizing the cost
\begin{eqnarray}
J(\Omega,  \mu_{\rm i} ) = {1\over 2} \sum_{k=1}^K \sum_{m=1}^M 
|u_{\Omega,  \mu_{\rm i}}(r_k,0,t_m) - d_k^m|^2,
\label{dcost}
\end{eqnarray}
where $u_{\Omega,  \mu_{\rm i} }(r_k,0,t_m)$ denotes the corresponding 
solution of the forward problem  evaluated at the receivers 
at the recording  times. More refined cost functionals  based on optimal transport 
\cite{yunan,metivier} could be employed. To reduce the occurrence of unphysical 
minima, the cost (\ref{dcost}) is regularized adding additional terms, terms of 
Tikhonov type, for instance, see \cite{jcp19, hohage2d} and section 
\ref{sec:linearized}.

To proceed, we need a mathematical representation for the geometry of the
anomalies in terms of a set of parameters $\nnu$.
Star-shaped parametrizations furnish a simple choice to represent the shape
of anomalies, though more general representations can be considered too
\cite{level_set, fichtner_object, palafox}. Star-shaped objects are defined
by a center and a radius function that fixes the position of the boundary 
points along all possible rays emerging from the center. 
Assume we know the tissue contains $L$ star-shaped anomalies,
that is, $\Omega = \cup_{\ell=1}^L \Omega^\ell$. Assume $\mu_{\rm i}(x)$ 
is piecewise constant, equal to $\mu_{\rm i}^\ell$ in $\Omega^{\ell}.$ 
Then, different representations of the radius function lead to lower or higher 
dimensional approaches. We will consider two possibilities.

For smooth star-shaped objects we can approximate the radius
function by trigonometric polynomials.
Given the data $\dn = (d_1^1,\ldots,d_K^1, \ldots, d_1^m,\ldots,d_K^m)$ 
we wish to predict  the $n(L,Q)=L(2Q+4)$ parameters
\begin{eqnarray} 
\nnu = (\nnu^1,\ldots,\nnu^L), \;\;
\nnu^\ell = (c_{x}^\ell, c_{y}^\ell,a_{0}^\ell,b_{1}^\ell,a_{1}^\ell,\ldots,
b_{Q}^\ell,a_{Q}^\ell,\mu_{\rm i}^\ell),  \;\;  \ell=1,\ldots,L, 
\label{parameters}
\end{eqnarray}  
representing the centers $(c_x^\ell,c_y^\ell)$ and radii $r^\ell(\theta)$ of 
the anomalies, ordered by blocks, associated to the parameterization 
\begin{eqnarray}
\label{trigonometric1}
\mathbf q(\theta)^\ell=
(c_x^\ell,c_y^\ell)+r^\ell(\theta)(\cos(2 \pi \theta),
\sin(2 \pi \theta)),\quad \theta \in[0, 1],\\
\label{trigonometric2}
r^\ell(\theta)= a_0^\ell + 2 \sum_{q=1}^Q a_q^\ell \cos(2 \pi q \theta) +
 2 \sum_{q=1}^Q b_q^\ell \sin(2 \pi q \theta),
\end{eqnarray}
for $\ell=1,\ldots,L$. Analogous parametrizations are available in three 
dimensions replacing Fourier expansions for the radius by expansions 
in terms of spherical harmonics  \cite{jcp19,hohage3d}. 
The number of modes $Q$ controls the allowed boundary roughness,
large values generate more complex shapes \cite{hohage3d}.

\begin{figure}[!hbt]
\centering
\includegraphics[width=10cm]{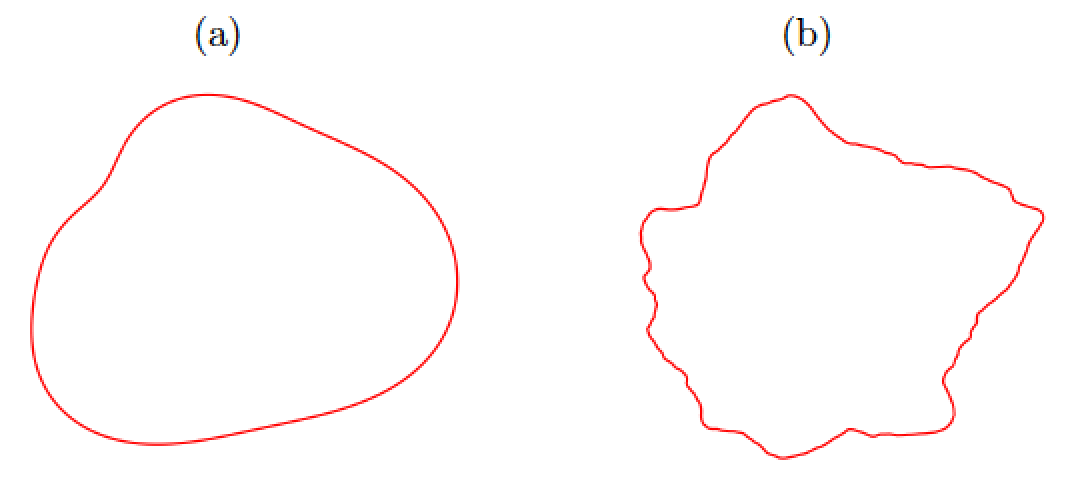} 
\caption{Star-shaped objects with radius defined by 
(a) a trigonometric polynomial (\ref{trigonometric2}) with $Q=5$ and 
(b) a piecewise approximation (\ref{trigonometric2_Z}) built from a uniform 
mesh $(\theta_j,r_j)$,  $j=0,\ldots, Z$, $Z=500$, with step $1/Z$.
}
\label{radius}
\end{figure}

Irregular boundaries are better represented by general radius functions 
$r(\theta)$, see \cite{babak, ghattas_radius}. In our case, we  approximate
the boundary by a piecewise linear reconstruction built on a uniform mesh
$\theta_j$ of $[0, 1]$ with node values $r_j=r(\theta_j)$, $j=0,\ldots, Z$.
The set of anomalies is then represented by the $n(L,Z)=L(3+Z)$
dimensional parameter set
\begin{eqnarray} 
\nnu = (\nnu^1,\ldots,\nnu^L), \quad
\nnu^\ell = (c_{x}^\ell, c_{y}^\ell,r_0^\ell,r_1^\ell,\ldots,
r_{Z-1}^\ell,\mu_{\rm i}^\ell),  \quad \ell=1,\ldots,L,
\label{parameters_Z}
\end{eqnarray} 
where $r_j^\ell=r^\ell(\theta_j)$, $j=0,\ldots, Z-1$,  and $r_Z=r_0$.
The boundary of each object is given by
\begin{eqnarray}
\label{trigonometric1_Z}
\mathbf q(\theta)^\ell=
(c_x^\ell,c_y^\ell)+r^\ell(\theta)(\cos(2 \pi \theta),
\sin(2 \pi \theta)),\quad \theta \in[0, 1],\\
\label{trigonometric2_Z}
r^\ell(\theta)=  r_j^\ell {\theta - \theta_{j+1} \over \theta_j - \theta_{j+1}}
+ r_{j+1}^\ell {\theta - \theta_{j} \over \theta_{j+1} - \theta_{j}}, \; 
\theta \in [\theta_j, \theta_{j+1}], \; j=0,1,\ldots,Z\!-\!1,
\end{eqnarray}
for $\ell=1,\ldots,L$. 
Notice that, while $Q$ is usually small, $Z$ can be very large.
Figure \ref{radius} compares a star-shaped object defined by 
(\ref{parameters})-(\ref{trigonometric2}), with a star-shaped object
defined by a piecewise approximation built from 
(\ref{parameters_Z})-(\ref{trigonometric2_Z}).

To quantify uncertainty in the solution of the inverse problem,
we resort to Bayes' formula \cite{KaipioSomersalo06,Tarantola05}
in finite dimension:
\begin{eqnarray}
p_{\rm pt}(\nnu):= p(\nnu | \dn) = { p( \dn | 
\nnu) \over p(\dn)} p_{\rm pr}(\nnu).
\label{bayes}
\end{eqnarray}
Here,  the prior density of the variables $p_{\rm pr}(\cdot)$  
incorporates available expert knowledge, while
$p(\dn|\nnu)$ represents the conditional probability (or likelihood)
of the observations $\dn$ given the variables $\nnu$.
The solution of the Bayesian inverse problem is the posterior density 
$p_{\rm pt}(\nnu|\dn)$  of the parameters given the data. The density 
$p(\dn)$ is a normalization factor that does not depend on the 
parameters. We choose a likelihood $p( \dn | \nnu)$ 
\begin{eqnarray}
p( \dn | \nnu) = {1 \over (2\pi)^{N/2} 
\sqrt{|\Gn|}} \exp \Big(- {1 \over 2} \|
\mathbf f( \nnu) - \dn  \|^2_{\Gn^{-1}} \Big).
\label{likelihood}
\end{eqnarray}
Here, $\| \mathbf v \|_{\Gn^{-1}}^2 =  
\mathbf {\overline v}^t \Gn^{-1} \mathbf v$ and 
$\mathbf f(\nnu)$ represents the measurement operator
associated to parameters $\nnu$, that is, 
\begin{eqnarray}
\mathbf f(\nnu)
=\left( U_{\Omega_{\nnu},\mu_{\rm i}}(\mathbf r_k,0,t_m) 
\right)_{k=1, \ldots, K, m=1, \ldots ,M},
\label{measurement}
\end{eqnarray} 
where  $U_{\Omega_{\nnu},\mu_{\rm i}}$ is the solution of 
the forward problem and $N=KM$.

We typically choose $p_{\rm pr}(\nnu)$ as a multivariate Gaussian
or a log Gaussian, see Section \ref{sec:prior} for details.
We could implement this approach using prior information obtained
by any means, for instance, other imaging systems or other imaging
algorithms, see \cite{transient}.  In the absence of this information, 
the next section explains how to generate prior knowledge 
from the data.

\section{Topological priors for the anomalies}
\label{sec:pgeometry}

Topological energy methods provide prior information on the number, location and
size of the anomalies by splitting the recorded data in two halfs: dodd represents the
data mesured at times $t_{2m+1}$ and deven the data measured at times $t_{2m}$: We exploit deven to generate prior information on the anomalies using the topological energy of the deterministic cost functional (\ref{dcost}). The remaining half of the data dodd enters the likelihood (\ref{likelihood}), as we will explain later.

\subsection{Calculation of topological energies}
\label{sec:te}

Given data $\deven$, the topological energy \cite{tenergy0, tenergy1} for 
the  cost
\begin{eqnarray}
J(\Omega) = {1\over 2} \int_{\Gamma_{obs}} \int_0^{\tau_{\rm end}}
| u_{\Omega}(\mathbf x,s) - \deven(\mathbf x,s) |^2 ds dx,
\label{tcost}
\end{eqnarray}
$u_\Omega$ being the solution of (\ref{forward_dim})-(\ref{bcneuman})
is given by
\begin{eqnarray}
E(\mathbf x) = \int_0^{\tau_{\rm end}} 
|U(\mathbf x,s) |^2 |P(\mathbf x,s)|^2 ds, 
\label{te}
\end{eqnarray}
where  $U= u_\Omega$ and $P$ is the associated adjoint field
that appears in the calculation of topological derivatives \cite{guzina2}.
In our set-up, we consider the observation set $\Gamma_{\rm obs}$
to be a set of receivers. Thus, $\int_{\Gamma_{\rm obs}}$ in (\ref{tcost}) 
becomes a sum of values at the receivers $\mathbf r_k$. Since we 
record data at discrete time values $t_{2m}$, we approximate 
$\int_0^{\tau_{\rm end}}$ by a sum  of values at such times too. 
Setting $\Omega= \emptyset$, the 
forward $U$ and adjoint  $P$ fields  are given by
\begin{eqnarray}
\begin{array}{ll}
 U_{tt} - c^2  \Delta U = f (t)  g(\mathbf x),  & \mathbf x \in R, \\[1ex]
{\partial U \over \partial \mathbf n} = 0,  & \mathbf x \in  \partial R, \\ [1ex]
U(\mathbf x, 0) = 0, U_t(\mathbf x, 0) = 0, & \mathbf x \in R, 
\end{array} \label{forward}
\end{eqnarray}
for $t\in [0, \tau_{\rm end}]$ and
\begin{eqnarray}
\begin{array}{ll}
[P_{tt} - c^2  \Delta P](\tau_{\rm end}\!-\!t) = - (U\!-\!\deven)(\tau_{\rm end}\!-\!t) 
\sum_{k=1}^K \delta_{\mathbf r_k},   & \mathbf x \in R, \\ [1ex]
{\partial P \over \partial \mathbf n} = 0,  & \mathbf x \in  \partial R, \\ [1ex]
P(\mathbf x, \tau_{\rm end}) = 0, P_t(\mathbf x, \tau_{\rm end}) = 0, & \mathbf x \in R,
\end{array}  \label{adjoint}
\end{eqnarray}
for $ t\in [\tau_{\rm end},0]$.
Here, $c$ is a constant equal to the healthy tissue wavespeed everywhere 
and $\delta_{\mathbf r_k}$ represent Dirac masses supported at the
receivers. For computational purposes, we replace them by Gaussian
regularizations.
Notice  that problems  (\ref{forward}) and (\ref{adjoint}) can be 
solved computationally even when $c_{\rm i}$ and $\mu_{\rm i}$ are unknown. 
This is an advantage over alternative methods based on topological 
derivatives \cite{guzina2} which require the knowledge of these parameters. 
The fact that spurious oscillations in the presence of multiple objects are 
considerably reduced constitutes an additional asset.
Topological energies are somewhat related to backpropagation techniques
\cite{papanico} and have been exploited for nondestructive testing of 
materials and tissues  in \cite{tenergy0,tenergy1,tenergybio}. 

\begin{figure}[!hbt]
\centering
\includegraphics[width=12cm]{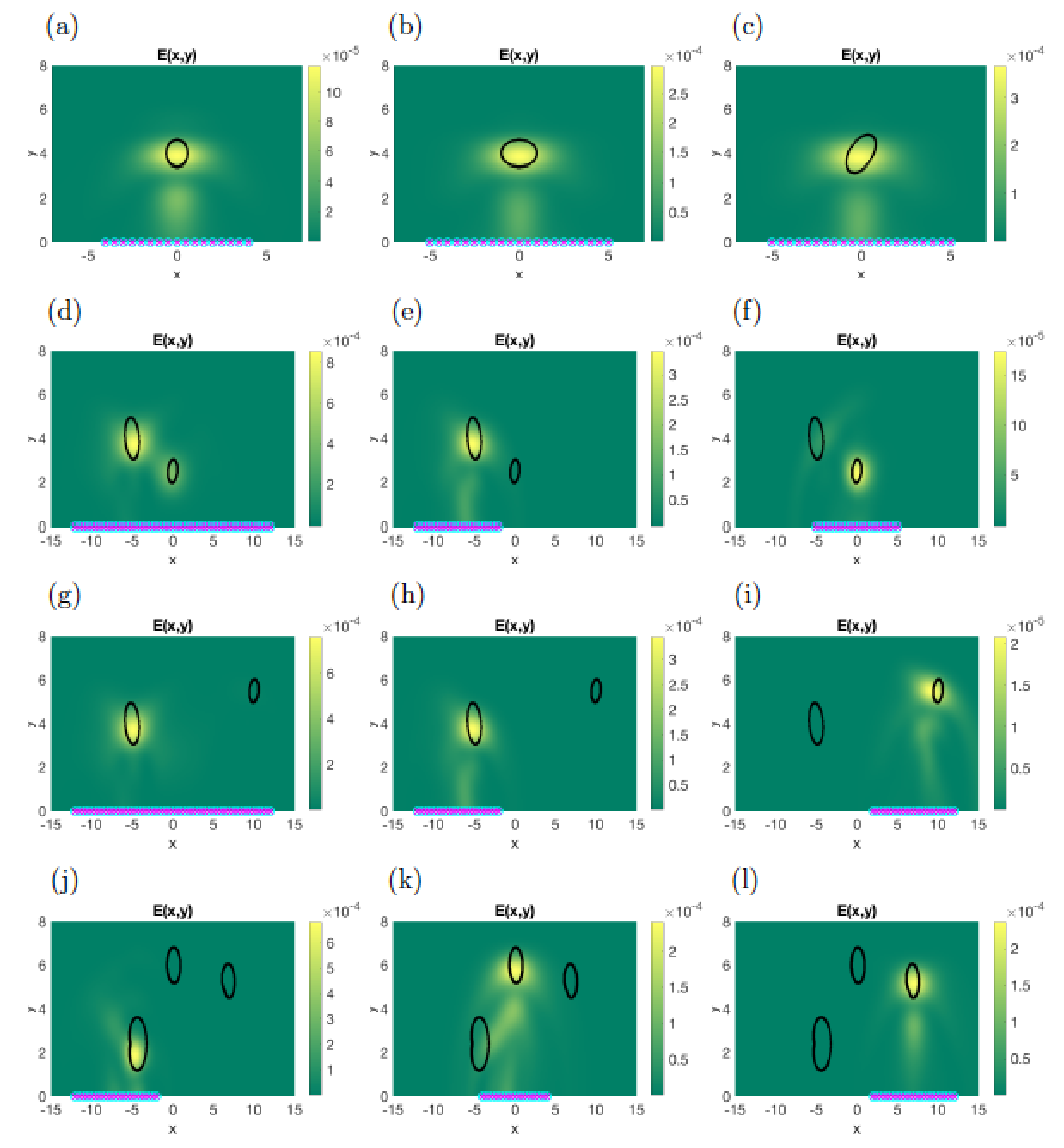}
\caption{Topological energy fields for: 
single objects with different sizes and orientations 
(a) circle, (b) ellipse, (c) rotated ellipse;
two objects under different emitter/receiver configurations
(d),(g) centered, (e),(h) left sided, (f),(i) right sided;
and three star-shaped objects (j)-(l) sweeping the bottom
region. 
Crosses and circles represent emitters and receivers, located 
at the same position. Black curves represent the true objects.
Noise level in the data: $10$ \%.}
\label{fig2}
\end{figure}

The previous description assumes that we record data at the
receivers from the  time $t=0$ at which we start to emit. 
If we start the recording later, at a time $\tau_{\rm in}$, formula 
(\ref{tcost})  integrates from $\tau_{\rm in}$ to $\tau_{\rm end}$ 
and the right hand side in (\ref{adjoint}) is only non zero in 
$[\tau_{\rm in},\tau_{\rm end}]$. We set the final time
$\tau_{\rm end} \sim 2 {H \over c}$, where $H$ is the expected 
resolution depth. 

Figure \ref{fig2} shows the topological energy fields obtained for 
several object geometries under different emitter/receiver configurations
for the parameter values specified in \ref{sec:nondimensional},
after removing dimensions.
The data $\deven$ used to calculate them are synthetic: they are 
generated by solving numerically the nondimensionalized forward 
problem (\ref{forward_num}) in the presence of the true objects, 
evaluating the solution in the selected space/time datagrid 
and adding $10$\% noise, as explained in Section 
\ref{sec:set-up}.
To prevent  inverse crimes, the fields $U$ in (\ref{forward}) and 
$P$ in (\ref{adjoint}) are approximated numerically using  rougher 
meshes: the spatial and time steps for them are twice the steps
used when solving numerically to generate the data, and the spatial
meshes vary.
We have set $\tau_{\rm in}=2$ (value at which $f(t)$ almost vanishes 
for our parameter choice) and $H=7$ in the calculation of the
topological energy.  

\subsection{Prior construction}
\label{sec:prior}

Once the topological fields are calculated, we construct
a first guess $\Omega_0$ for the anomalies immersed in a background 
medium  $R$ by exploiting the peaks of the topological energy:
\begin{eqnarray}
\Omega_0 = \left\{ \mathbf x \in R \; \vert \; E(\mathbf x) >
(1-C_0) \, {\rm max}_{\mathbf y \in R} E(\mathbf y)
\right\},
\label{omega0}
\end{eqnarray}
where  $C_0 \in (0,1)$  is such  that $J(\Omega_0) < J(\emptyset).$
In case several objects are present, we obtain more precise information by
sequentially activating fractions of the whole network of emitters/receivers,
as shown in Figure \ref{fig2}(d)-(l). We fit circles to the dominant peaks found 
for each fraction. In this way, we are able to detect all the anomalies. Instead, 
when we use the information coming from the whole network, we often find 
the most  prominent anomaly only. 
We use this information to construct priors for the two types
of star-shaped parameterizations we consider as follows.

Assuming we locate $L$ peaks, we fit  to them circles parametrized by
$\nnu_0 = (\nnu^1_0,\ldots,\nnu^L_0)$. When we work with the representation
(\ref{parameters})-(\ref{trigonometric2}), we set
\begin{eqnarray}
\nnu^\ell_0 = (c_{x,0}^\ell, c_{y,0}^\ell,a_{0,0}^\ell,b_{1,0}^\ell,a_{1,0}^\ell,\ldots,
b_{Q,0}^\ell,a_{Q,0}^\ell,\mu_{\rm i,0}^\ell),  \quad \ell=1,\ldots,L, 
\label{nuprior}
\end{eqnarray}
where $c_{x,0}^\ell, c_{y,0}^\ell$ is the center of mass of each component,
$a_{0,0}^\ell$ half the smallest diameter, and
$a_{1,0}^\ell=\ldots=a_{Q,0}^\ell=b_{1,0}^\ell=\ldots=b_{Q,0}^\ell=0$. We 
also set $\mu_{\rm i,0}^\ell = \mu$, the known background value for the
healthy tissue. Then, we choose $p_{\rm pr}(\nnu)$ as a multivariate 
Gaussian with covariance matrix $\Gpr$ 
\begin{eqnarray} \begin{array}{l}
p_{\rm pr}(\nnu)  
={1\over (2 \pi)^{n/2}}  {1\over \sqrt{|\Gpr |}}
\exp(-{1\over 2}
(\nnu - \nnu_0)^t \Gpr^{-1} 
(\nnu - \nnu_0) ), 
\end{array} \label{prior}
\end{eqnarray}
where $n$ is the dimension of $\nnu_0$,
provided that $\mu_{\rm i}^\ell>0$, the curves associated 
to the parameterization $\nnu$ fulfill $r^\ell(\theta)>0$, for 
$\theta \in [0,1]$ and all $\ell$, they do not intersect, and they 
do not form nested configurations. Otherwise, $p_{\rm pr}(\nnu)$
is set equal to zero. Notice that $r^\ell(\theta)>0$ is not a condition 
on the sign of the curve parameters, but on the sign of the 
combination (\ref{trigonometric2}).
We choose a  diagonal covariance matrix $\Gpr$ formed by  $L$ 
blocks. In our numerical tests, each block starts with 
$(\sigma_x^\ell)^2=(\sigma_y^\ell)^2=0.1$
and ends with $(\sigma_\mu)^2=20^2.$
Then $(\sigma_{a_0}^\ell)^2=0.1$ and $(\sigma_{a_q}^\ell)^2 
= (\sigma_{b_q}^\ell)^2 = 0.1/(1+q^2)^s$, 
$ 1 \leq q \leq Q$, $s$ large, as in \cite{sergei}, so that the prior favors 
regular shapes with $r(t)>0$. Typically, we fix $s=3$ and $Q=5$.

When we work with the representation 
(\ref{parameters_Z})-(\ref{trigonometric2_Z}), we set
\begin{eqnarray}
\nnu^\ell_0 = (c_{x,0}^\ell, c_{y,0}^\ell,r_{0,0}^\ell,r_{1,0}^\ell,\ldots,
r_{Z-1,0}^\ell,\mu_{\rm i,0}^\ell),  \quad \ell=1,\ldots,L, 
\label{nuprior_Z}
\end{eqnarray}
where $c_{x,0}^\ell, c_{y,0}^\ell$ is the center of mass of each 
component, $r_{0,0}^\ell = r_{1,0}^\ell = \ldots = r_{Z-1,0}^\ell$ is
half the smallest diameter and $\mu_{\rm i,0}^\ell = \mu$. 
Notice that $r_{Z,0}=r_{0,0}$.
We choose $p_{\rm pr}(\nnu)$ as the product of  multivariate
Gaussians for the variables $c_{x,0}^\ell$, $c_{y,0}^\ell$,
$\mu_{\rm i,0}^\ell$, with the same standard deviations as
before, and log Gaussian distributions for $r_{0,0}^\ell,r_{1,0}^\ell, 
\ldots, r_{Z-1,0}^\ell,$ with Matern covariance matrices, see  
\cite{matt_thesis}. We use for the Matern covariance between 
points $(\cos(\theta_i),\sin(\theta_i))$ and
$(\cos(\theta_j),\sin(\theta_j))$ separated by a distance $d$ the 
expressions \cite{matern}
\begin{eqnarray*}
C_{\nu,\rho,\sigma}(d) = \sigma^2 {2^{1-\nu} \over \Gamma(\nnu)} 
\left( \sqrt{2\nu}{d\over \rho} \right)^\nu K_\nu \left( \sqrt{2\nu}{d\over \rho}\right),
\label{matern}
\end{eqnarray*}
where $\Gamma$ is the Gamma function, $K_\nu$ the modified
Bessel function of the second kind, and $\nu$, $\rho$, $\sigma$ are
parameters.  In our numerical tests we fix $\nu=3/2$, $\sigma=0.2$,
$\rho=0.5$.

\section{Markov Chain Monte Carlo sampling}
\label{sec:MCMC}

We insert the prior distributions obtained in the previous 
section in the posterior probability $p_{\rm pt}$ given by (\ref{bayes}) 
and  (\ref{likelihood}), with the data $\dodd$ not used to produce the
prior information. Then we can sample the unnormalized posterior 
distribution $q(\nnu) =  p( \dodd | \nnu)  p_{\rm pr}(\nnu)$  using Markov 
Chain Monte Carlo (MCMC) methods. Note that the  unknown scaling  
factor $p(\dn)$ in (\ref{bayes}) is not needed for MCMC sampling. 
Standard MCMC methods, such as Metropolis-Hastings or Hamiltonian
MonteCarlo \cite{mcmc}, produce a chain of $N$-dimensional states 
$\nnu^{(0)} \longrightarrow  \nnu^{(1)} \ldots \longrightarrow  \nnu^{(i)} 
\ldots$  which evolve to be distributed according to the target distribution.
One first samples an initial state $\nnu^{(0)}$ from the prior 
distribution,  and then moves from one state $ \nnu^{(i)}$ to 
the next $ \nnu^{(i+1)}$ guided by a transition operator.  
More recent  ensemble MCMC samplers \cite{matt_sampler,goodmanweare} 
draw $W$ initial states (the `walkers' or `particles') from the prior distribution 
and transition to new states while mixing them to construct the chain. This 
allows us to handle  multimodal posteriors \cite{sergei}  and to parallelize 
the process for faster  exploration of the structure of the posterior distribution. 

Different ensemble samplers adapt better to the different parametrizations
we consider for the anomalies. Affine-invariant samplers perform well
in our set-up. We have considered two.  
The first one is a  stretch move based Affine Invariant Ensemble Sampler 
(SAIES), see \cite{goodmanweare}:
\begin{itemize}
\item Initialization: Choose $W$ initial states $\boldsymbol \nu_w^{(1)}
\in \mathbb R^d$, $w=1, \ldots, W,$ with probability $\pi$ (the prior 
probability $p_{\rm pr}$ in our case) and a value $a>1$.  
\item For each step $s=1,\ldots,S$, 
\begin{itemize}
 \item For each $w=1, \ldots, W$
 \begin{itemize}
   \item Draw $\boldsymbol \nu_q^{(s)}$ at random from the set 
   $\{ \boldsymbol \nu_j^{(s)} \}_{j \neq w}$.
   \item Choose a random value $z_w$ from the distribution $g(z)= {
   1\over \sqrt{z}}$  when $z \in [1/a,a]$, zero otherwise.
   \item Set $\boldsymbol \nu_{w,\rm prop}^{(s)} =
   \boldsymbol \nu_{w}^{(s)} + z_w
  (\boldsymbol \nu_w^{(s)} - \boldsymbol \nu_q^{(s)}).$
   \item  Set 
   $\boldsymbol \nu_{w}^{(s+1)}=\boldsymbol \nu_{w,\rm prop}^{(s)}$
   with probability ${\rm min}\left\{1, z_w^{d-1} { 
    p_{\rm pt}(\boldsymbol \nu_{w,\rm prop}^{(s)})
    \over  p_{\rm pt}(\boldsymbol \nu_{w}^{(s)}) } \right\},$
   or else keep $\boldsymbol \nu_{w}^{(s+1)}=\boldsymbol \nu_{w}^{(s)}.$
 \end{itemize}
\end{itemize}
\item Output: The samples  $\boldsymbol \nu_w^{(s)}$, $w=1, \ldots, W$,
$s=1,\ldots,S$.
\end{itemize}
The second one is a general Affine Invariant Ensemble Sampler (AIES), which 
proceeds as follows, see  \cite{matt_sampler} for instance:
\begin{itemize}
\item Initialization: Choose $W$ initial states $\boldsymbol \nu_w^{(1)} 
\in \mathbb R^d$, $w=1, \ldots, W$, with probability $\pi$ (the prior probability 
$p_{\rm pr}$ 
in our case) and a value $\lambda>0$.
\item For each step $s=1,\ldots,S$, 
\begin{itemize}
 \item For each $w=1, \ldots, W$
 \begin{itemize}
   \item Set $\overline{\boldsymbol \nu}  = {1\over W-1} \sum_{j\neq w}  
   \boldsymbol \nu_j^{(s)}$.
   \item Draw $z_w$ with probability ${\cal N}(0,1).$
   \item Set $\boldsymbol \nu_{w,\rm prop}^{(s)} =
   \boldsymbol \nu_{w}^{(s)} + {\lambda \over \sqrt{W-1}}
   \sum_{j\neq w} z_j (\boldsymbol \nu_j^{(s)} - \overline {\boldsymbol \nu}).$
   \item  Set 
   $\boldsymbol \nu_{w}^{(s+1)}=\boldsymbol \nu_{w,\rm prop}^{(s)}$
   with probability ${\rm min}\left\{1, { p_{\rm pt}(\boldsymbol 
   \nu_{w,\rm prop}^{(s)}) \over p_{\rm pt}(\boldsymbol \nu_{w}^{(s)})} \right\},$
   or else keep $\boldsymbol \nu_{w}^{(s+1)}=\boldsymbol \nu_{w}^{(s)}.$
 \end{itemize}
\end{itemize}
\item Output: The samples  $\boldsymbol \nu_w^{(s)}$, $w=1, \ldots, W$,
$s=1,\ldots,S$.
\end{itemize}

\begin{figure}[!hbt]
\centering
\includegraphics[width=10cm]{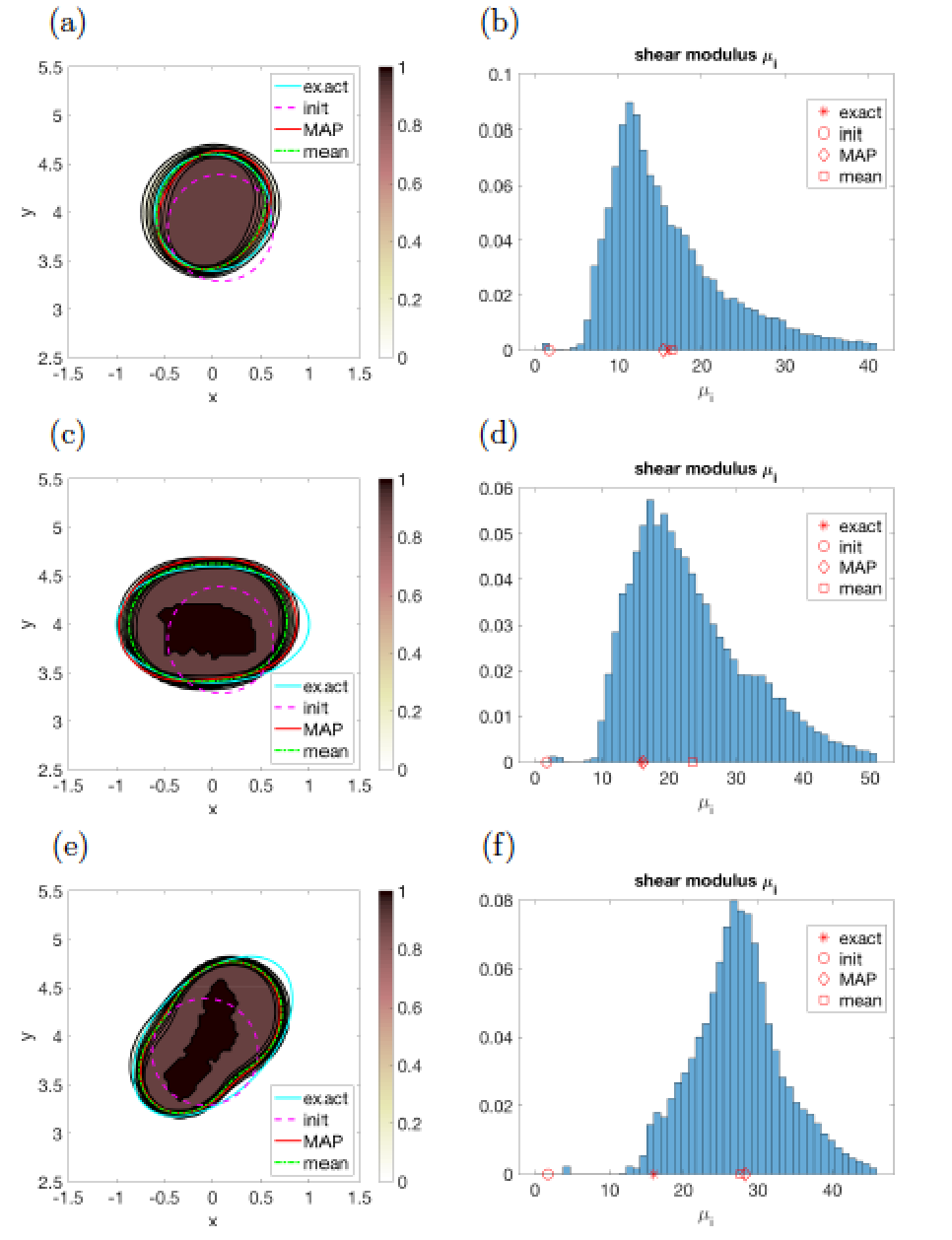} 
\caption{(a), (c), (e) True objects versus MAP estimate and sample mean 
calculated from MCMC samples for different geometries. The contour levels 
represent the probability of belonging to the object.
(b), (d), (f)  Histograms quantifying uncertainty of the MAP estimate and mean
values for $\mu_{\rm i}$.
Parameters and samplers: SAIES with $a=2$, $\tilde S=500$,
$W=480$ and $B= W \tilde S /5$.}
\label{fig3}
\end{figure}

\begin{figure}[!hbt]
\centering
\includegraphics[width=12cm]{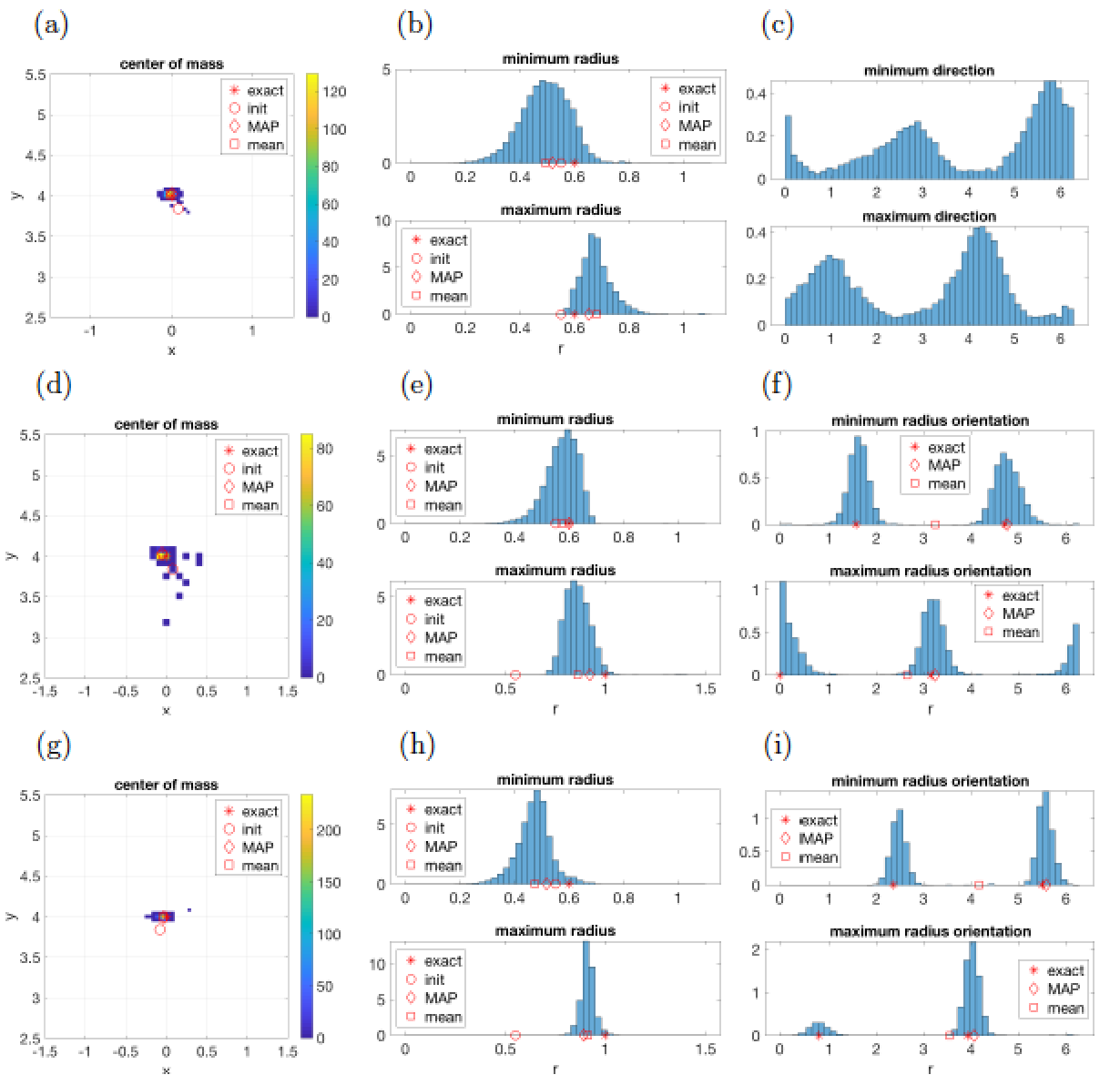}
\caption{Histograms representing a discrete approximations of the densities 
for the distribution of the centers of mass (a), (d), (g),  radius size
(b), (e), (h), and orientation (c), (f), (i) in the three test geometries
considered in Figure \ref{fig3}.
Same sampling parameters.
}
\label{fig4}
\end{figure}

\begin{figure}[!hbt]
\centering
\includegraphics[width=12cm]{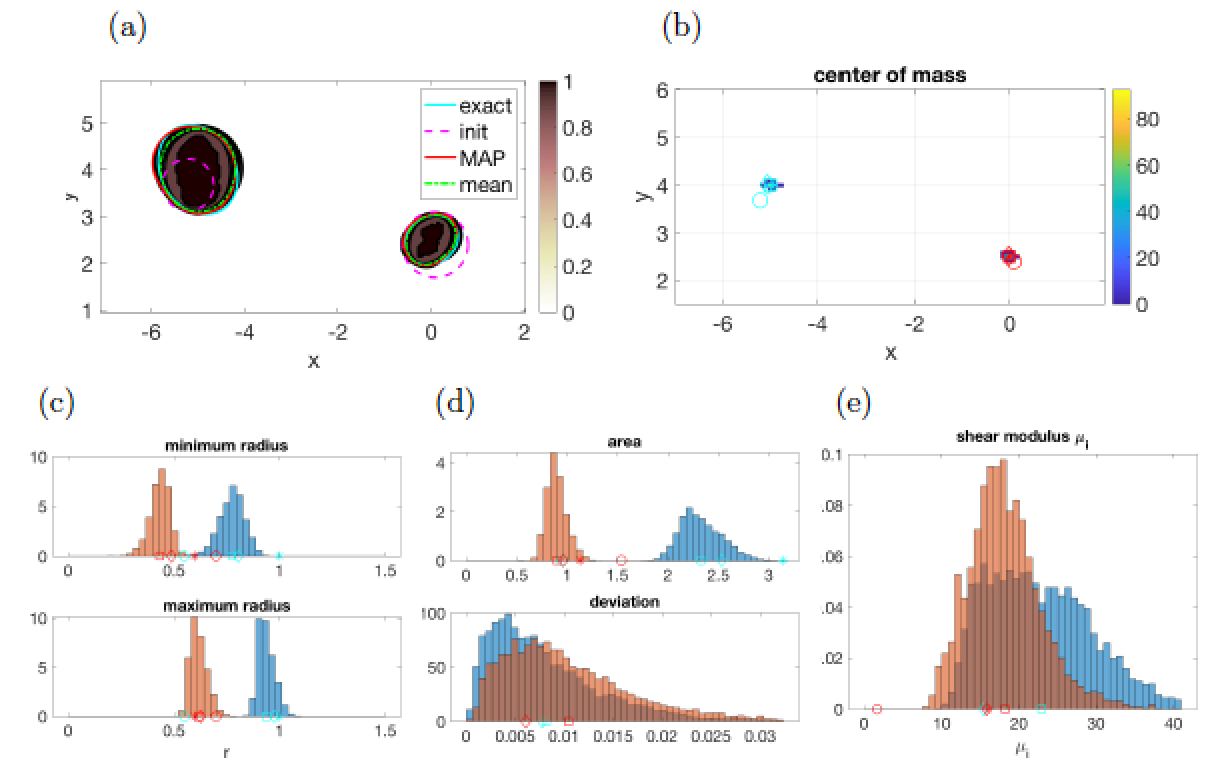}
\caption{
(a) True objects versus MAP estimate and mean calculated from MCMC samples. 
The contour levels represent the probability of belonging to the object.
Superimposed curves represent the exact contours (solid light cyan), the 
MAP point (solid dark red), the mean (dash-dotted green), and the initial guess 
(dashed magenta).
(b) Histograms representing discrete approximations of the densities 
for the distribution of the centers of mass.
(c)-(d) Histograms for the radius sizes (c), area and deviation from a circular 
object (d) and shear modulus $\mu_{\rm i}$ (e).
Blue histograms and cyan symbols correspond to the large object, 
orange histograms and red symbols to the small one
(asterisk: exact value, circle: initial value, diamond: MAP point, square: mean).
Parameters and samplers: SAIES with $a=2$, $\tilde S=500$,  $W=480$ 
and $B= W \tilde S /5$. }
\label{fig5}
\end{figure}

\begin{figure}[!hbt]
\centering
\includegraphics[width=12cm]{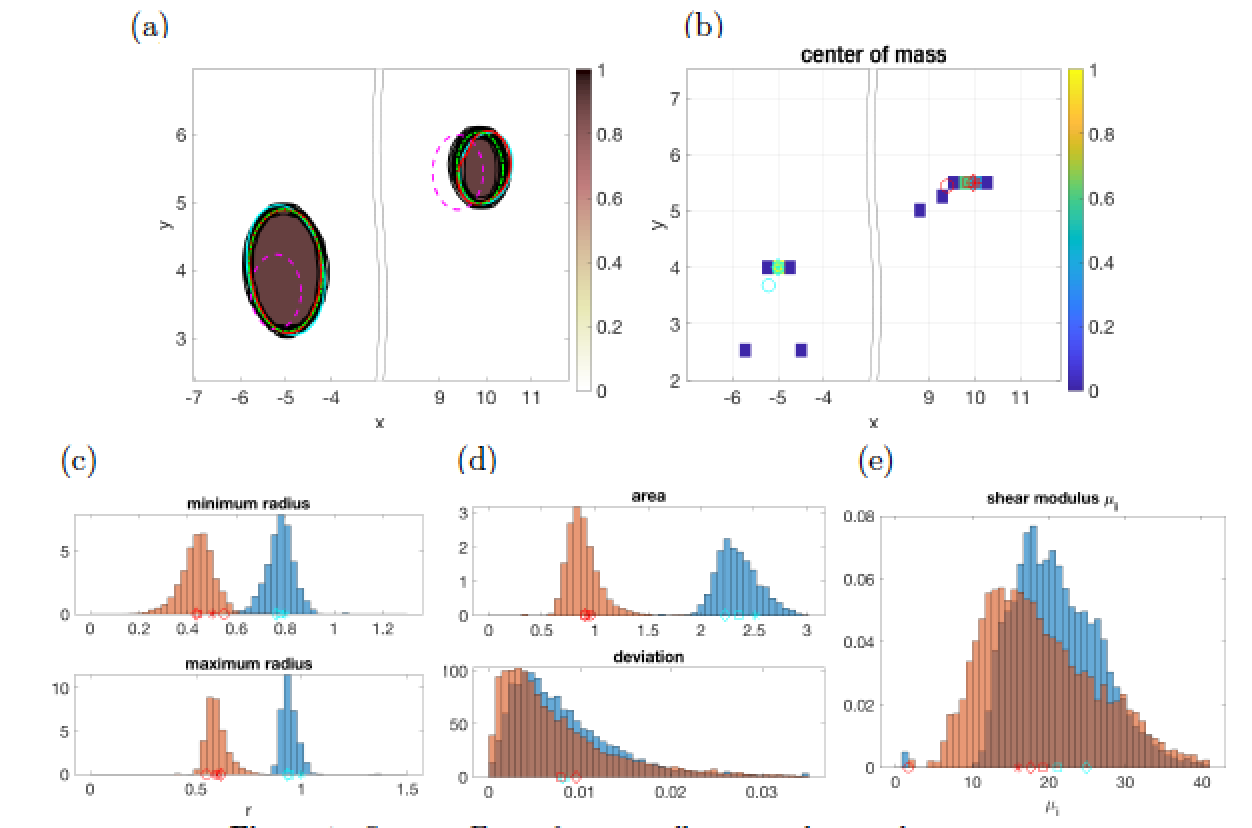} 
\caption{Same as Fig. \ref{fig5} for two well separated anomalies.}
\label{fig6}
\end{figure}

 \begin{figure}[!hbt]
\centering
\includegraphics[width=12cm]{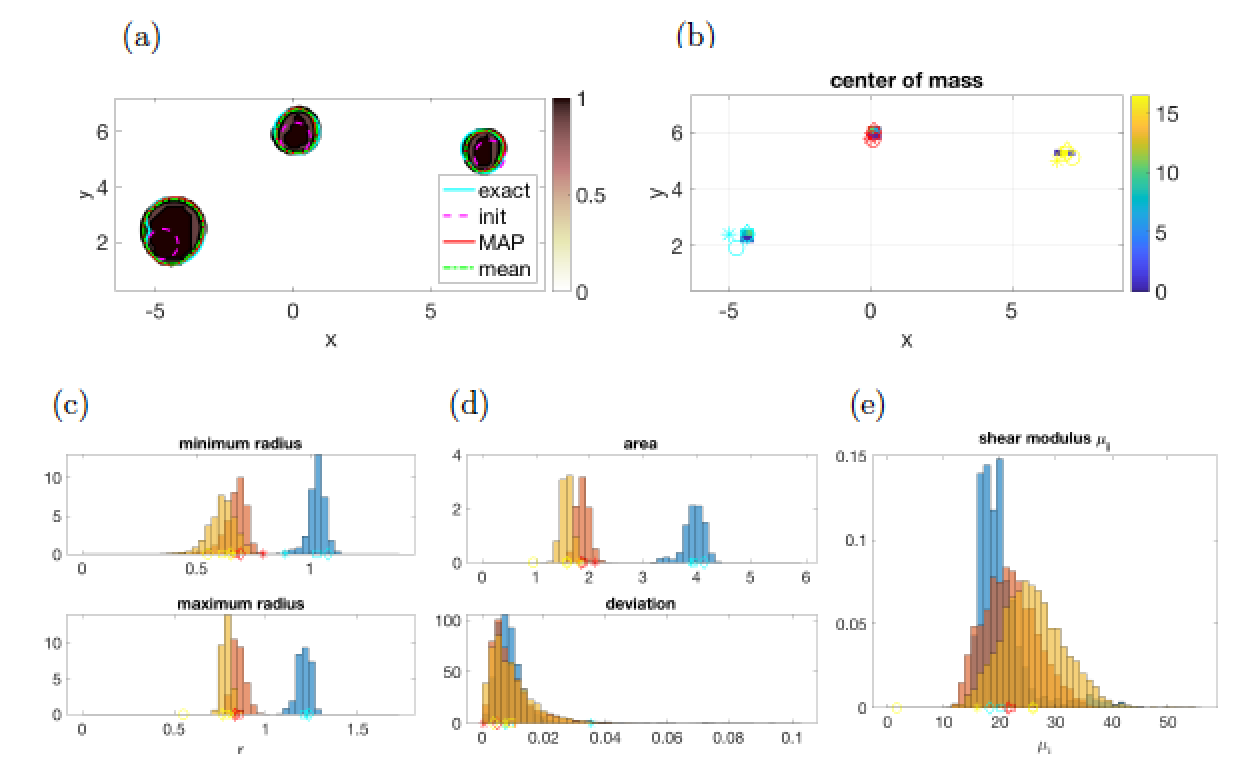}
\caption{Same as Figure \ref{fig5} for a configuration with three objects.
Orange histograms and red symbols correspond to the middle object, 
blue histograms and cyan symbols to the left-most one, 
yellow histograms and  symbols to the right-most one.}
\label{fig7}
\end{figure}

While the first sampler evolves faster in low dimensions $d$, it usually requires 
$W>2d$ to perform properly.  In principle, the general  AIES  can be more 
robust as dimension grows. 
Figures  \ref{fig3}-\ref{fig7} display results with SAIES for smooth shapes 
admitting low dimensional parametrizations. AIES provides similar
results doubling the number of steps. Section \ref{sec:irregular} considers 
high dimensional irregular shapes. There, AIES perfoms reasonably well 
with $W$ slightly larger than $d$. 
In both cases and for each $w$, we keep one of each three samples up to 
a total number of $\tilde S = S/3$ to reduce correlations and discard  the 
first $\tilde S /5$ as a  burn in period.
We have set $\sigma_{\rm noise}= \alpha
 \, {\rm max} |d^k_{n, {\rm true}}| /100$ in $\Gn$ with $\alpha =10$.

Here, for highly smooth shapes, we consider the parameter set (\ref{nuprior}) and
define $\Gpr$ as  in Section \ref{sec:prior}. Then, we insert the prior probability 
(\ref{prior}) obtained by topological methods and the likelihood (\ref{likelihood}) 
in the posterior probability (\ref{bayes}) to be sampled.
From the samples, we obtain information on the most likely values 
for $\nnu$, that is, the maximum a posteriori (MAP)  estimate and
the uncertainty about it, depicted in figures
\ref{fig3}-\ref{fig7}. 
Figure \ref{fig3} illustrates the uncertainty in the shape of the
anomaly and the value of the parameter representing the 
dimensionless shear modulus  $\mu_{\rm i}$ for single shapes:
a circle and an ellipse with different orientations. The sample with 
highest probability defines the MAP point
and the mean of the parameters corresponding to all the samples 
defines a mean estimate. The location and shape of the anomalies is
reasonably well captured by both, see also Figure \ref{fig4} for the 
uncertainty in geometrical features of interest, such as the location 
of the center of mass, the size of the largest and smallest diameters 
and their orientation. However, the value of the shear modulus
displays larger uncertainty, still in the range indicating sickness. 
The histograms reveal distributions with wide and asymmetric tails.
Notice that a change in the orientation of an object can drastically 
increase uncertainty in the predictions, compare Fig. \ref{fig3}(d) and
Fig. \ref{fig3}(f). 

Figures \ref{fig5}-\ref{fig7} consider configurations with multiple
anomalies. The approximation of the 
shapes provided by the MAP  point and the mean values in figures 
\ref{fig5}-\ref{fig7} is quite reasonable, regarding both the shapes,
their basic geometrical features and the shear modulus, though
we observe again wide asymmetric tails.

When we include in the prior less anomalies than needed,
the distribution may be multimodal: we may spot the missing
components.
When we include in the prior more anomalies than needed,
the spurious ones may essentially vanish because $\mu_{\rm i}$
is basically  equal to $\mu$. 
Notice that our priors contained the correct number of anomalies.
The resulting distributions represent a single mode. Also, 
the means and the MAP points are reasonably close. This suggests
that optimization schemes could capture the MAP estimate, allowing 
for a Laplace approximation of the posterior distribution, which
can be sampled at a much lower cost. The computational time
drops from a few days to a  few minutes.

\section{Sampling from a Bayesian linearized formulation}
\label{sec:linearized}

To reduce the computational cost, we analyze the Laplace 
approximation of the posterior density (\ref{bayes}) obtained by 
linearization at the maximum a posteriori (MAP) point $\nnumap$.
This strategy first computes the vector $\nnumap$, 
which minimizes the negative log likelihood 
\begin{eqnarray}
J(\nnu) = {1\over 2 \sigma_{\rm noise}^2} \sum_{m=1}^M
\sum_{k=1}^K |U_{\boldsymbol  \nu}(r_k,0,t_m) \!-\! d_{k}^m|^2
+ {1 \over 2} (\boldsymbol  \nu - \nnu_{0})^t
\Gpr^{-1} (\boldsymbol  \nu - \nnu_{0}),  
\label{cost_nu}
\end{eqnarray}
where $U_{\nnu}$ is the solution of the forward problem with  
object $\Omega$ parametrized by $\nnu$
given by (\ref{parameters}). Then, we approximate the posterior
distribution by a multivariate Gaussian ${\cal N}(\nnumap,\Gpo)$ with 
posterior convariance matrix $\Gpo = \mathbf H_{\nnumap}^{-1}$, where 
$\mathbf H_{\nnumap}$ is an approximation of the Hessian of the 
measurement operator (\ref{measurement}) evaluated at $\nnumap$
\cite{georg_linearized,georg_marmousi}.

\subsection{Computing the MAP point}
\label{sec:cMAP}

The MAP point is calculated exploiting techniques of deterministic
optimization. Taking the prior $\nnu_{0}$ as initial guess of the
parametrization, that is, $\nnu^0 = \nnu_{0}$, we can implement the 
Newton  type  iteration $ \nnu^{j+1} = \nnu^{j} + 
\boldsymbol \xi^{j+1}$ where  $ \boldsymbol \xi^{j+1}$ is the solution of  
\begin{eqnarray}
\left(\mathbf H(\nnu^j) + \omega_j
{\rm diag}(\mathbf H(\nnu^j) ) \right)
\boldsymbol \xi^{j+1} = -  \mathbf g(\nnu^j),
\label{it_false}
\end{eqnarray}
see \cite{lmf}, where $\mathbf H(\nnu)$ and $\mathbf g(\nnu)$
represent the Hessian and the gradient of the cost.
In practice, to reduce the occurrence of negative radii and the risk of 
loop formation in the curves, we introduce an additional parameter 
$\lambda>0$, replacing (\ref{cost_nu}) by
\begin{eqnarray}
J_{\lambda}(\nnu) = {1\over 2 \sigma_{\rm noise}^2} \sum_{m=1}^M
\sum_{k=1}^K |U_{\boldsymbol  \nu}(r_k,0,t_m) \!-\! d_{k}^m|^2
\!+\! {\lambda \over 2} (\boldsymbol  \nu - \nnu_{0})^t
\Gpr^{-1} (\boldsymbol  \nu - \nnu_{0})
\label{cost_lambda}
\end{eqnarray}
and (\ref{it_false}) by
\begin{eqnarray}
\left(\mathbf H^{\rm GN}_{\lambda_j}(\nnu^j) + \omega_j 
{\rm diag}(\mathbf H^{\rm GN}_{\lambda_j}(\nnu^j) ) \right)
\boldsymbol \xi^{j+1} = - \mathbf g_{\lambda_j}(\nnu^j).
\label{it_true}
\end{eqnarray}
Here,  the subscript $\lambda_j$  indicates that we multiply $\Gpr^{-1}$
by a factor $\lambda_j$ in the initial iterations to balance the two terms 
defining the cost $J$ in (\ref{cost_nu}). Notice that we have also
replaced the full Hessian by the Gauss-Newton part of the Hessian
to reduce the computational cost per iteration. The components
of $\mathbf H^{\rm GN}_{\lambda}(\nnu)$ and 
$\mathbf g_{\lambda}(\nnu)$ are given by:
\begin{eqnarray}
(g_{\lambda}(\nnu))_i & =  {\partial J(\nnu) \over \partial \nu_i} & =  
{1\over \sigma^2_{\rm noise}} \sum_{m=1}^M \sum_{k=1}^K 
(U_{\nnu}(r_k,0,t_m) - d_k^m) 
{\partial U_{\nnu} \over \partial  \nu_i}(r_k,0,t_m)  \nonumber \\
 & & +  \lambda [\Gpr^{-1}(\nnu - \nnu_0)]_i, \label{grad} \\
(\mathbf H^{\rm GN}_{\lambda}(\nnu))_{i,\ell} & =  
{\partial^2 J(\nnu) \over \partial  \nu_i \partial \nu_\ell} & \sim   
{1\over \sigma^2_{\rm noise}}
\sum_{m=1}^M \sum_{k=1}^K {\partial U_{\nnu}\over 
\partial \nu_i}(r_k,0,t_m)  {\partial U_{\nnu} \over \partial \nu_\ell}(r_k,0,t_m)
\nonumber \\
& &  +\lambda  [ \Gpr^{-1}]_{i,\ell}. \label{hess}
\end{eqnarray}
The second order derivatives of $U_{\nnu}$ are neglected. 

To optimize our objective function  we implement a double 
iteration:
\begin{itemize}
\item Initially, we set $\omega_0=10^{-4}/2$, $\lambda_0= 0.1 
\sigma_{\rm noise}^{-2}$, and $\nnu^0 = \nnu_0$.
\item At each step we  calculate $\nnu^{j+1} = \nnu^{j} + \boldsymbol \xi^{j+1}$, 
where  $\boldsymbol \xi^{j+1}$ is the solution of (\ref{it_true}). Then
\begin{itemize}
\item We check
i) if $r(\theta)>0$, $r(\theta)$ given by (\ref{trigonometric2}), and if
$\mu_{\rm i} > 0.5 \mu$,
ii) if the  functional $J_{\lambda_j}(\nnu)$ decreases replacing 
$ \nnu^{j}$ with $\nnu^{j+1}$.
\item If any of these conditions fails, we do not accept $\boldsymbol \xi^{j+1}$. 
We increase $\omega_j$ by a factor  $2$, solve again (\ref{it_true}) and
check conditions i) and ii) until they are fulfilled.
\item If both conditions are satisfied, we accept $\boldsymbol \xi^{j+1}$ 
and set  $\omega_{j+1}=\omega_j/2$ and $\lambda_{j+1} = 
{\rm max}(\lambda_j/5,1)$.
\end{itemize}
\item After a few steps $j_0$, $\lambda_{j+1}=1$ for $j \geq j_0$.
When the relative difference between the new value of 
the  cost and the previous one is smaller than a tolerance Tol (here
Tol $=5 \times 10^{-7}$), we freeze all the components except
$\mu_{\rm i}^{j+1}$ and iterate with respect to $\mu_{\rm i}$ until
variations fall below a threshold $0.02$.
\end{itemize}

To evaluate the derivatives ${\partial U_{\nnu} \over \partial \nu_i}(r_k,0,t_m)$
required for the calculation of (\ref{grad})-(\ref{hess}) at each step,
we use the approximation
\[
{\partial U_{\nnu} \over \partial \nu_i}(r_k,0,t_m) 
\sim {U_{\nnu+\eta_i}(r_k,0,t_m) - 
U_{\nnu}(r_k,0,t_m)  \over \eta_i}
\]
with $\eta_i$ small,  $U_{\nnu+\eta_i}$ being the  solution of the forward
problem with $\nu_i$ replaced by $\nu_i+\eta_i.$  All the forward problems
are solved with the same discretization and steps we used in Section
\ref{sec:MCMC}.
The values of $\eta_i$ must be calibrated.
Initially, we set for each block $\ell=1,\ldots, L$ in (\ref{parameters})
$\eta_1^\ell=\eta_2^\ell=\eta_3^\ell=\eta_{2Q+4}^\ell=\eta = 0.1$ and
$\eta_{3+2i}^\ell=\eta_{2+2i}^\ell= \eta/2$ for $i=1,\ldots,Q$.
As we iterate, we calibrate values for $\eta_i$ estimating
the quotients ${ D_{\nu_i} U_{\nnu^j}(r_k,0,t_m)
\over D^2_{\nu_i} U_{\nnu^j} (r_k,0,t_m)} $, where
$D$ and $D^2$ represent approximations of derivatives,
and averaging over $k$ and $m$. 
In the tests we have performed, the choice
\begin{eqnarray*}
\eta_1^\ell=\eta_2^\ell = 0.05, \quad \eta_{2Q+4} = 0.15, \\
\eta_{3+2i}^\ell=\eta_{2+2i}^\ell= 0.05, \quad i=1,\ldots,Q,
\label{eta}
\end{eqnarray*}
gives good results, with $\eta_3$ in the range $0.025-0.225$. 

\begin{figure}[!hbt]
\centering
\includegraphics[width=12cm]{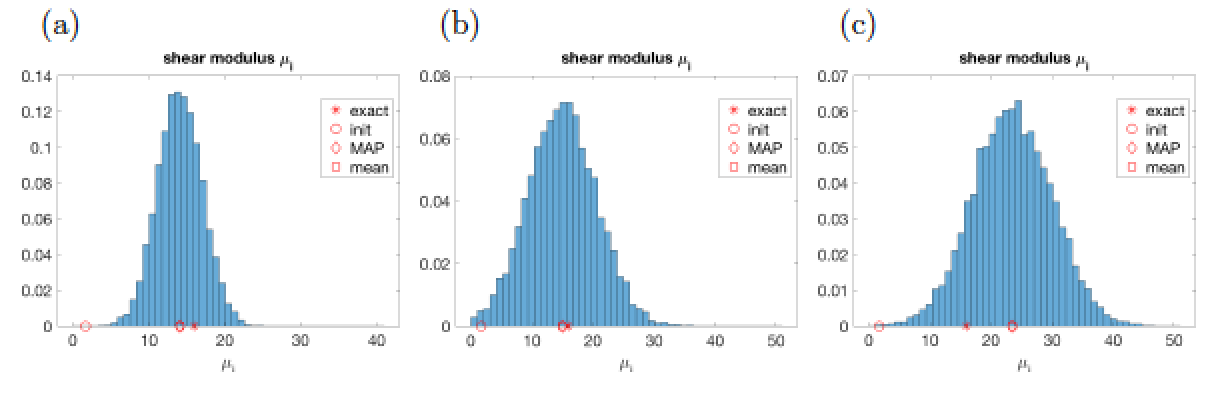}
\caption{Counterpart of Fig. \ref{fig3} (b), (d), (f) obtained by calculating the
MAP point and linearizing the posterior probability about it. 
10000 samples plotted.
}
\label{fig9}  
\end{figure}

\begin{figure}[!hbt]
\centering
\includegraphics[width=9cm]{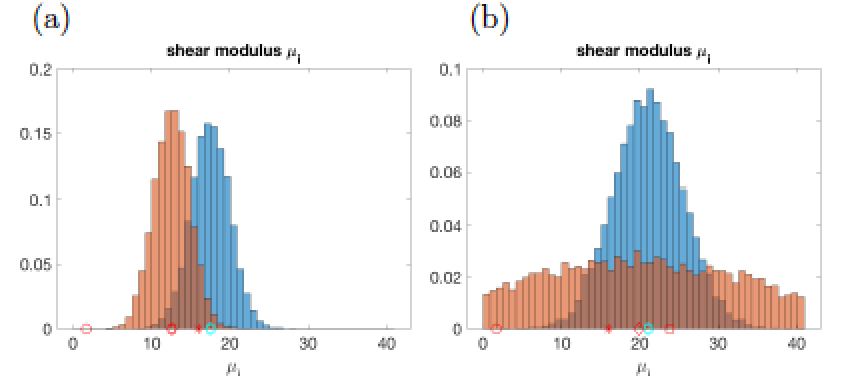}
\caption{Counterparts of Fig. \ref{fig5}(b) and Fig. \ref{fig6}(b) obtained by linearized Bayesian
methods. $10000$ samples plotted. 
}
\label{fig1112}
\end{figure}

\begin{figure}[!hbt]
\centering
\includegraphics[width=12cm]{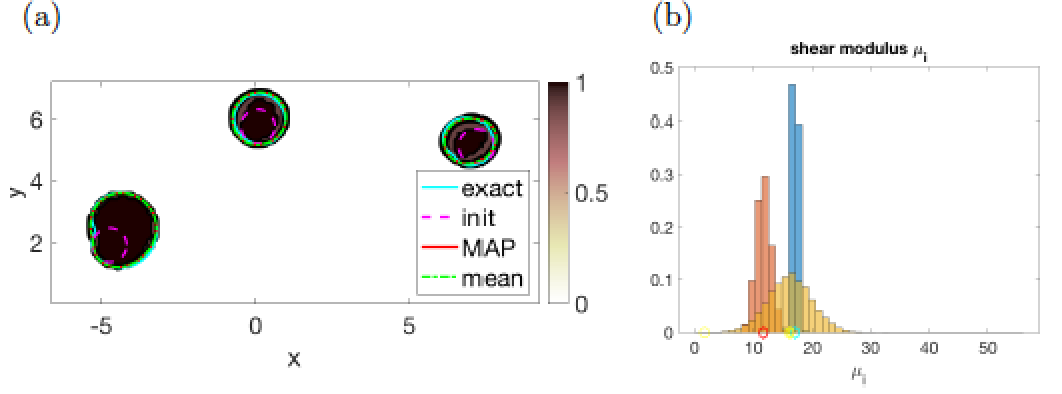}
\caption{Counterpart of Fig. \ref{fig7} obtained by linearized Bayesian
methods. $10000$ samples plotted. 
}
\label{fig13}
\end{figure}

\subsection{Sampling}
\label{sec:bsampling}

Once we have obtained an approximation to $\nnumap$,
we linearize the posterior distribution about it, approximate
by a multivariate Gaussian distribution ${\cal N}(\nnumap, \Gpo)$
and  draw samples from it to quantify uncertainty. 
We set
\begin{eqnarray*}
\Gpo = ( \mathbf F(\nnumap)^t \Gn^{-1} 
\mathbf F(\nnumap) + \Gpr^{-1})^{-1} = \mathbf H^{\rm GN}(\nnumap)^{-1},
\end{eqnarray*}
where $\mathbf F(\nnumap)= \left({\partial U_{\nnumap} \over \partial 
\nu_i}(p_j) \right)_{j,i}$ 
and 
\[ \mathbf p = ((r_1,0,t_1), \ldots, (r_K,0,t_1), \ldots,
(r_1,0,t_M), \ldots, (r_K,0,t_M)), 
\]
that is, $\mathbf F(\nnumap)$ is the matrix with $i$th-column
${\partial U_{\boldsymbol  \nu}\over \partial \nu_i}(r_k,0,t_m) $,
$k=1,\ldots,K$, $m=1,\ldots,M$ evaluated at $\nnumap$.
This can be done by means of the relation
\begin{eqnarray}
\nnu = \boldsymbol  \nu_{\rm MAP} + \Gpo^{1/2} 
\mathbf n, \label{sample_map}
\end{eqnarray}
$\mathbf n$ being a standard normal randomly generated vector (iid).

Figures \ref{fig9}-\ref{fig13} revisit the previous MCMC tests with this 
procedure.
In each case, we optimize to approximate the MAP point and generate a 
large collection  of samples of the posterior distribution by means 
of (\ref{sample_map}). The values of the cost for the approximated MAP
estimates obtained this way are similar to those for the MAP estimates 
previously found by MCMC sampling. 
Comparing the results, we remark that the MAP points and contour curves
for the shapes remain similar. Fig \ref{fig13}(a) illustrates this fact in the
example with three objects. However, the values of $\mu_{\rm i}$ show 
larger variability. For single objects, the MAP points and mean values
remain alike, while the wide asymmetric tails are lost. The tests with
more objects show a similar tendency. Notice that as objects become
smaller and distant, information can be lost, as it happens in the red
histogram in Fig \ref{fig1112}(b) (compare to Fig \ref{fig6}(b)).

The computational cost of this approach is much smaller. The MAP estimate for
single objects is obtained in about $10$ steps, about $5$ minutes in a laptop
using MATLAB, while MCMC sampling can take 2-4 days depending on the
size of the computational regions.

\section{Irregular shapes}
\label{sec:irregular}

Finally, we consider irregular shapes defined by high dimensional
parametrizations of the form (\ref{parameters_Z})-(\ref{trigonometric2_Z}). 
We insert the prior distributions (\ref{nuprior_Z}) obtained by
topological methods in the posterior probability given by (\ref{bayes}) 
and (\ref{likelihood}), with the data $\dodd$ not used to produce the
prior information.

\begin{figure}[!hbt]
\centering
\includegraphics[width=10cm]{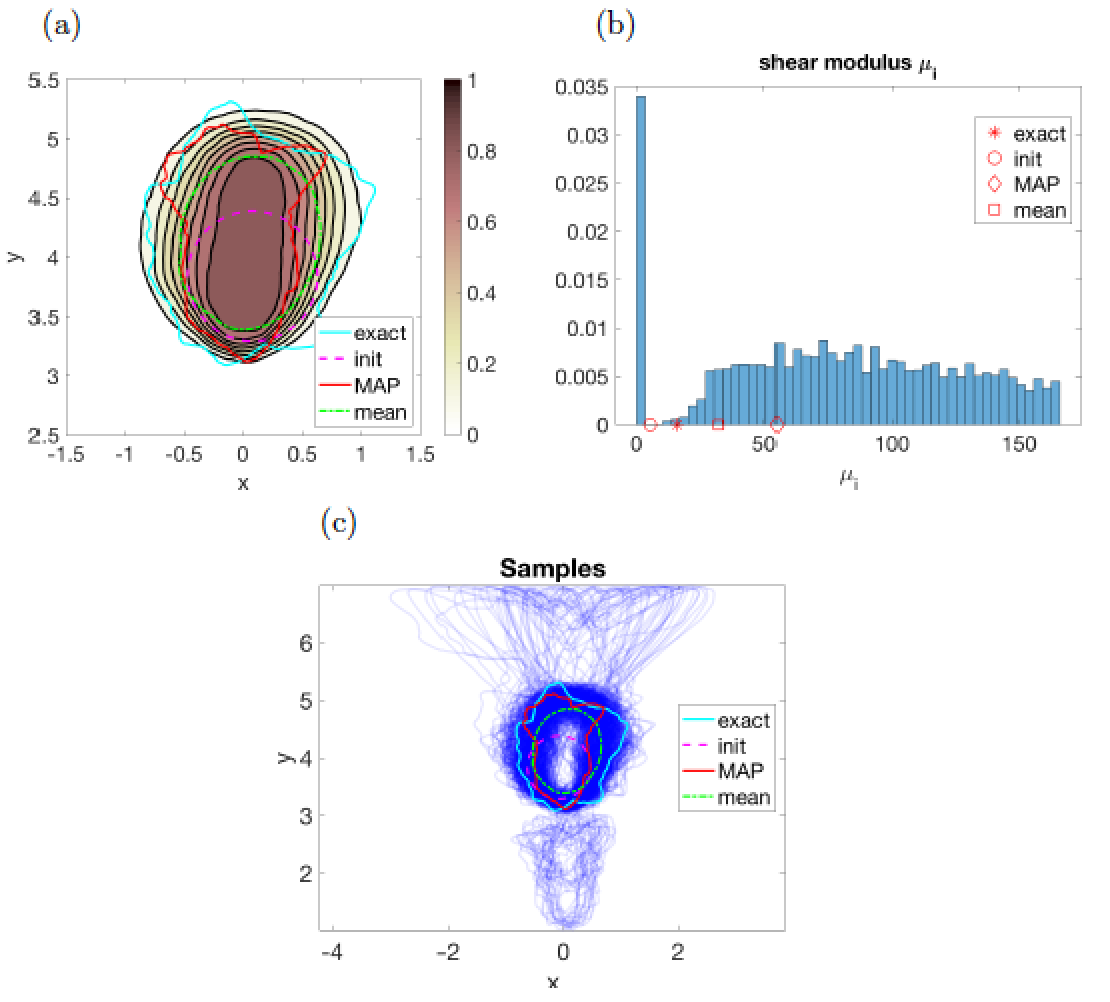} 
\caption{Results for the object in Figure \ref{fig2}(b) for $Z=500$.
MCMC with AIES $W=600$, $\tilde S=3900$, $B= W \tilde S/2$ , $\lambda=0.2$.}
\label{fig15}
\end{figure}

The affine-invariant ensemble sampler AIES described in Section 
\ref{sec:MCMC} produces the results represented in Figure \ref{fig15} 
for an irregular shape. Notice that the use of SAIES would require 
$W > 2d = 1006$ walkers, which would mix much more slowly. 
Now, the MAP estimate does not approach the true shape. Nevertheless, 
the mean profile and the contour plot give an idea of the location and 
size of the anomaly.

In the previous sections, the probability for negative $\mu_{\rm i}$ was 
set equal to zero. Now, $\mu_{\rm i} = \exp(\gamma)$, where $\gamma$ 
is the random variable that we sample. Notice the peak for $\mu_{\rm i}$ 
near zero. It is due to a family of large samples with small $\mu_{\rm i}$.
Figure \ref{fig15}(c) represents the last $W$ samples we obtained.
We observe a dominant family of samples which wrap around the true 
object. A second family is formed by smaller shapes with larger values
of  $\mu_{\rm i}$  placed between the object and the emitters, at the 
location of a small secondary peak of the topological energy 
(see Figure \ref{fig2}(c)).  The third family corresponds to large samples
with small $\mu_{\rm i}$ placed behind the true object.
The sample distributions we obtain this way are multimodal, though the
main mode dominates the rest when averaging to obtain a mean.
Notice that uncertainty in the values of $\mu_{\rm i}$ with this procedure
seems quite large.
The prior information we use has low quality in this case. We look
for an irregular shape assuming that the prior is a smooth circle
and $\mu_{\rm i}$ is the value for the healthy tissue. Inconsistency
between the prior and the data may lead to multimodality, as pointed
out in the previous section.

 \begin{figure}[!hbt]
\centering
\includegraphics[width=10cm]{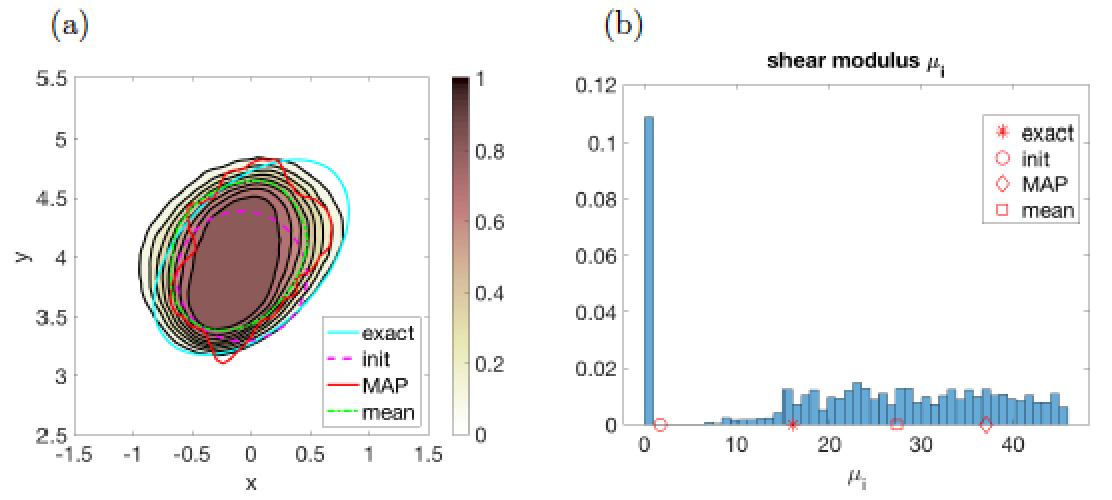}
\caption{Results for the rotated ellipse, working with high dimensional
parametrizations allowing for irregular shapes with $Z=500$.
Same sampling parameters as in Fig. \ref{fig16} except $\tilde S=2400$.
}
\label{fig16}
\end{figure}

Working with smooth shapes, we get Figure \ref{fig16} for the
rotated ellipse already studied in Figure \ref{fig3}(c). The samples
we generate behave in a similar way as those in Figure \ref{fig15}(c).
The MAP point is unlikely to be smooth when we do not enforce
the prior knowledge we have on smoothness.  However, the information 
provided by the  mean parameters and the statistics of geometrical 
characteristics and values for shear moduli is still useful, though less 
precise.
Enforcing a smooth parametrization we get better results for smooth 
shapes, at a lower computational cost, see
Figure \ref{fig3}(c). Similarly, the irregular shape studied in Figure
\ref{fig15} could be studied in the smooth framework employed
in  Section \ref{sec:MCMC} to obtain information about mean values
at a lower cost.

\section{Conclusions}
\label{sec:conclusions}
We have developed a Bayesian approach for the detection and characterization 
of anomalies in tissues which uses topological energies to generate priors. 
In this framework, anomalies are represented by star-shaped objects 
whose shear moduli differ from the surrounding tissues. We have 
considered low dimensional parametrizations for simple smooth shapes 
and higher dimensional approximations for irregular shapes, which
can be used to distinguish encapsulated (smooth) and invasive 
(irregular) tumors, for instance.

For simple shapes, MCMC methods based on different types of affine 
invariant ensemble samplers provide a good characterization of the 
structure of the posterior distribution, which displays asymmetric tails 
for each mode representing an object. This approach is time consuming, 
since we must generate a few hundred thousand samples by
solving a time dependent wave equation for each of them.
We have shown that its is possible to approximate the `maximum a 
posteriori'  (MAP)  estimate of the parameters defining the hardness 
and geometry of these anomalies. To do so, we minimize a proper cost 
functional, which can be done in a few iterations by Newton type
iterations. Linearizing the parameter-to-observable map about the MAP 
point, we are able to quantify the uncertainty in nature of the anomalies, their 
location and shape by  generating samples of the Laplace approximation 
to  the posterior distribution at a low  computational cost.

We have tested these schemes in 2D shear imaging set-ups, finding 
reasonable agreement between both sampling techniques for such
shapes.
While MCMC sampling furnishes a deeper insight in the structure 
of the posterior, including asymmetry and possible multimodality, the 
linearization approach provides results quite fast. This is essential 
for potential technological applications and three dimensional extensions.
However, it may miss multimodality and asymmetry details.

Irregular shapes lead to higher dimensional problems and optimization
approaches to calculate a MAP point encounter difficulties due to 
fast variations in the boundary. We have seen that affine invariant samplers 
which are robust as dimension grows still provide some information
on basic anomaly properties, though we identify multimodality
features due to inconsistency between the data and the prior.
Better descriptions of the anomaly shape and shear modulus would 
probably require improved prior knowledge or a different type of 
parametrization.

Alternative Bayesian formulations seek variations in the wave speed of 
the whole tissue, which leads to infinite dimensional problems and very 
large computational cost.  The approach based on seeking  shapes
described by a moderate number of parameters that we propose here 
has been tried for simple shapes on time independent imaging problems 
for which efficient boundary element solvers are available. 
Lacking similar solvers for time dependent wave problems, we have
succeeded in developing fast finite element schemes allowing us to
implement our Bayesian formulation in terms of parametrized boundaries,
at a low computational cost, which is convenient for practical applications.

\appendix
\section{Approximate solutions for the forward problem}
\label{sec:discretization}

We recall here the pertinent existence and regularity result
for the forward problem, as well as some discretization
details and parameter choices.

\subsection{Existence and regularity}
\label{sec:ap_existence}

In the sequel, $H^1$, $H^2$ represent the standard Sobolev 
spaces and $(H^{1})'$ is the dual space of $H^1$ \cite{adams,brezis}.
$L^2$ stands for the usual space of square-integrable functions. \\

{\bf Theorem 1.} {\it Let $R$ and $\Omega$ be $C^1$ domains,  
$\Omega \subset R$ \footnote{The result remains
true with piecewise boundary regularity or when $R$ is a convex Lipschitz 
domain using Sobolev space theory for them \cite{adams,necas}.}. 
Assume $f \in C^\infty(\mathbb R^+) \cup L^\infty(\mathbb R^+)$ 
and $g \in C^\infty(\mathbb R^2) \cup L^\infty(\mathbb R^2)$.
Then, the problem (\ref{forward_dim})-(\ref{bcneuman}) has
a unique solution $u \in C([0,\tau];H^1(R))$, $u_t \in C([0,\tau];L^2(R))$, 
$u_{tt} \in L^2(0,\tau;(H^1(R))')$, for any  $\tau>0$.  
Furthermore, if $u_{tt}(\mathbf x, 0)=0$,  we also have 
$u_t \in C([0,\tau];H^1(R))$ and 
$u \in C([0,\tau];H^2(R\setminus \overline \Omega))$.} 

{\bf Proof.}
Existence of a solution $u$ with the stated regularity for wave equations
with positive and bounded $\mu$ and $\rho$ is a particular case of
results established in \cite{lions,raviart}. 
If $u_{tt}(\mathbf x, 0)=0$, $u_t$ solves (\ref{forward_dim})-(\ref{bcneuman}) 
with $f$ replaced by $f'$. Hence, $u_t \in C([0,\tau];H^1(R))$ and 
$u_{tt} \in C([0,\tau];L^2(R)).$
Then equation (\ref{forward_dim}) implies that $\Delta u(t) \in L^2(R
\setminus \overline \Omega)$, thus $u(t) \in H^2(R
\setminus \overline \Omega)$ by elliptic regularity theory
and $u$ is defined on 
$\Sigma$ and the receiving sites both in the sense of $L^2(\Sigma)$ 
traces and pointwise \cite{brezis,haraux}.

\subsection{Physical parameters and nondimensionalization}
\label{sec:nondimensional}

For computational purposes, we nondimensionalize the problem using
characteristic times and lengths. Let $T$ and $L$ be two characteristic
time and length scales to be chosen.
To simplify, one can take $\rho_{\rm i} = \rho \sim 1000 \, \rm kg/m^3 $
in tissues, though $\rho_{\rm i} > \rho$ in general (slightly).
We set $\mathbf x= \mathbf x' L$, $t= t' T$,  $u = u' L$, 
$\Omega = \Omega' L$, $R =  R' L$ and $\Sigma =  \Sigma' L$. 
Making the change of variables  and dropping the symbol $'$ for ease 
of notation, we get
\begin{eqnarray}
\begin{array}{ll}
u_{tt} - {\rm div} ({\mu T^2 \over \rho L^2}  \nabla  u) = {T^2 \over \rho L}
f(t T)  G(\mathbf x L) = \tilde f(t) \tilde G(\mathbf x), & \mathbf x \in R, \; t>0, \\[1ex]
 {\partial u \over \partial \mathbf n} = 0, & \mathbf x \in   \partial R, \\[1ex]
u(\mathbf x, 0) = 0, u_t(\mathbf x, 0) = 0, & \mathbf x \in R.
\end{array} \label{forward_adim}
\end{eqnarray}
Here, $\tilde f(t)= f_0 {T^2 \over \rho L}  (1- 2 \pi^2 f_M^2 T^2 t^2) 
e^{-\pi^2 f_M^2 T^2 t^2}$. 
We choose $\tilde G$  to have zero normal derivative at the interface,
for instance, $\tilde G(\mathbf x ) =  {1 \over ( \pi \kappa )^{n/2} } \sum_{j=1}^J
\exp(- {|(\mathbf x - \mathbf x_j)|^2 \over  \kappa })$, $n=2$. This function
represents the location of the emitters.

Typical experimental conditions \cite{tumorspeed, tumorelasticity, tumorstiffness}  
suggest the choice $L =1$ cm = $10^{-2}$ m and $T = 10^{-2}$ s.
For instance, typical anomaly shapes and sizes in a liver framework are ellipsoids  
of about $0.963 \times 1.15$ cm, buried at a depth between $6$ and $12$ cm. 
To spot anomalies of size $1$ cm, that is, $10^{-2}$ m, we should need a receiver 
grid of step about  $10^{-3}$ m distributed or moving over regions of cm length.
Typical parameter ranges \cite{tumorspeed, tumorelasticity, tumorstiffness} are 
$\mu_{\rm i} = 96 - 241$ kPa (carcinoma), $\mu_{\rm i} = 55 -71 $ kPa (normal 
tissue),  and $\mu_{\rm i} = 36 - 41 $ kPa (bening hyperplasia) in a prostate gland, 
for instance.
In a liver framework, $\mu_{\rm i} = 0.4 - 6$ kPa (healthy tissue) and $\mu_{\rm i}
= 15 - 100$ (unhealthy tissue). Breast is less appropriate for these methods because 
carcinoma may yield $\mu_{\rm i} = 22-560$ kPa, ovelapping with fibrous tissue 
$\mu_{\rm i} = 96-244$ kPa, normal fat $\mu_{\rm i} = 18-24$ kPa, and normal 
gland  $\mu_{\rm i} = 28-66$ kPa, other techniques \cite{fichtner_breast} may 
be more suitable.
Frequencies $f_M$ in shear elastography devices are $4-15$ Hz, or $50$ Hz, 
or $100-300$ Hz, depending on sizes involved.
                                                                                           
\begin{table}[!hbt]
\centering
\begin{tabular}{|c|c|c|c|c|c|c|c|c|c|}
\hline
$L$ & $T$ & $\rho_{\rm i}$ & $\rho$ & $\mu_{\rm i}$ & $\mu$ & $c_{\rm i}$ 
& $c$ & $f_M$ & $f_0$ \\
\hline
$0.01$ m & $0.01$ s & $\rho$ & $10^3 \rm {kg \over m^3}$ & $16$ kPa
& $1.69$ kPa  & $4 \rm {m\over s}$ & $1.3 {m\over s}$ & $50$ Hz & 
${\rho L \over T^2}$ \\
\hline
\end{tabular}
\caption{Dimensional parameters used in the simulations.}
\label{table:parameters}
\end{table}                                                                                          
                                                                                                                                                  
In our numerical tests we work with the parameters listed in 
Table \ref{table:parameters}.
We set $\mu_{\rm i} =16$ kPa and $\mu= 1.69$ kPa,
which results in wave speeds $c_{\rm i} = 4$ m/s inside the
anomalies and $c = 1.3$ m/s outside, a low contrast situation.
We select $f_0$ such that $f_0 {T^2 \over \rho L} =1$
and  $f_M=50$ Hz so that $f_M T = 0.5$.
Then, the final dimensionless forward problem is
\begin{eqnarray}
\begin{array}{ll}
u_{tt} - {\rm div} (c(\mathbf x)^2  \nabla  u)  =  
\tilde f(t) \tilde G(\mathbf x), 
& \mathbf x \in R,   \\ [1ex]
{\partial u \over \partial \mathbf n} = 0, & \mathbf x \in   \partial R, 
\\ [1ex]
u(\mathbf x, 0) = 0, u_t(\mathbf x, 0) = 0, & \mathbf x \in R,
\end{array} \label{forward_num}
\end{eqnarray}
for $t>0$, with
\begin{eqnarray*}
c^2(\mathbf x) = {\mu(\mathbf x) T^2 \over \rho L^2} = 
\left\{ \begin{array}{ll} 
1.69,    & \mathbf x \in R \setminus \overline{\Omega}, \\
16, & \mathbf x \in  \Omega,
\end{array}  \right.  \quad
c(\mathbf x) =  \left\{ \begin{array}{ll} 
1.3,    & \mathbf x \in R \setminus \overline{\Omega}, \\
4, & \mathbf x \in  \Omega,
\end{array}  \right.  
\label{coefs_adim} 
\end{eqnarray*}
and
\begin{eqnarray}
\tilde f(t) \tilde G(\mathbf x) = 
(1\!-\! 2 \pi^2 0.5^2 t^2) { e^{-\pi^2 0.5^2 t^2}
 \over \pi \kappa } \sum_{j=1}^J
e^{- {|\mathbf x - \mathbf x_j|^2 \over  \kappa }}. 
\label{source}
\end{eqnarray}

We generate synthetic data for our simulations by solving 
numerically (\ref{forward_num})-(\ref{source}) for different 
choices of anomalies $\Omega$ and adding random noise. 
We have used finite elements \cite{DautrayLions, raviart} with 
spatial step $\delta x=0.08$ and a total explicit spatial discretization 
with time  step $\delta t=0.00125$, see next section for  details.
We locate emitter/receivers at  fixed grids of step $0.5$ (or $0.2$)
and record the signal at a fixed time grid of step $0.025$.
The value of $\kappa$ can be adjusted to the step $\delta x$, so that  
it affects just a few nodes around the emitter. Here, we have set 
$\kappa =2$.
Alternatively, one could also perform an even extension at the interface  
$\Sigma = \{(x,y) \, | \, y=0 \}$ to get a problem set in the whole 
space and resort to boundary elements for wave problems 
\cite{tonatiuh} representing the emitters as point sources. However,
an adequate framework to implement such boundary value 
approach is still missing.

\subsection{Discretization}
\label{sec:ap_discrete}

To reduce the computational cost we focus on a limited tissue
region and truncate the computational region in such a way that 
$R$ is a rectangular region, as in Figure \ref{fig1}.  On the artificial 
boundaries $\partial R \setminus \Sigma$, we will enforce non 
reflecting boundary conditions \cite{nonreflecting}. On $\Sigma$, 
we keep the zero Neuman condition.
For the spatial discretization,  we use $P_1$ finite elements  on a 
fixed mesh of step $\delta x$ in space \cite{DautrayLions,raviart}. 
If $V= {\rm span}\{\phi_1,\ldots,
\phi_D\} \subset H^1$ is the resulting finite element space, we 
approximate $u$ by $u^D = \sum_{i=1}^D a_i(t) \phi_i $. Therefore, 
we must find $u^D$ such that
\begin{eqnarray*}
 \int_R u^D_{tt}(\mathbf x,t) \phi_j(\mathbf x) d \mathbf x  
 + \int_R c(\mathbf x)^2 \nabla u^D(\mathbf x,t) \phi_j(\mathbf x) d \mathbf x  
 - \int_{\partial R \setminus \Sigma} c^2 {\partial u^D\over \partial \mathbf n}
 \phi_j(\mathbf x) d S_{\mathbf x}  \\
 = \tilde f(t) \int R \tilde G(\mathbf x) \phi_j(\mathbf x) d \mathbf x,
 \end{eqnarray*}
for $j=1,\ldots,D$. Next, we use the nonreflecting boundary condition 
${\partial u^D\over \partial \mathbf n} \sim - {1\over c} u_t^D$ on  
$\partial R \setminus \Sigma$ and total discretizations for the time 
derivatives on a time mesh $t_n$ of step $\delta t$ \cite{DautrayLions,raviart}
\begin{eqnarray*}
u^D_{tt}(\mathbf x,t_n) \sim {u^D(\mathbf x,t_{n+1}) - 2u^D(\mathbf x,t_n)
+ u^D(\mathbf x,t_{n-1}) \over \delta t^2 }, \\
u^D_{t}(\mathbf x,t_n)  \sim {u^D(\mathbf x,t_n) - u^D(\mathbf x,t_{n-1}) 
\over \delta t}.
\end{eqnarray*}
To calculate the coefficients $a_i(t_n)$, $i=1,\ldots D$, we solve the recurrence 
relations 
\begin{eqnarray*}
\sum_{i=1}^D M_{j,i} a_i(t_{n+1}) =  
\sum_{i=1}^D  M_{j,i} (2 a_i(t_{n}) - a_i(t_{n-1})) 
-  \delta t^2  \sum_{i=1}^D  A_{j,i} a_i(t_{n})  \\
- c \, \delta t \sum_{i=1}^D B_{j,i}   (a_i(t_{n})- a_i(t_{n-1}))  
+ \delta t^2 \tilde f(t_n) G_j,
\end{eqnarray*}
for $n \geq 1$, where $M_{j,i} = \int_R \phi_j \phi_i d \mathbf x$,
$A_{j,i} = \int_R c^2 \nabla \phi_j \nabla \phi_i d \mathbf x$,
$B_{j,i} = \int_{\partial R \setminus \Sigma}  \phi_j \phi_i 
d S_{\mathbf x}$, $G_j = \int_R  \tilde G \phi_j d \mathbf x.$
The coefficients $a_i(t_0)=0$ and $a_i(t_1)=0$ for $i=1,\ldots,D$ 
are determined using the initial conditions.
The numerical solutions defined in this way are continuous, so that
the costs (\ref{dcost}), (\ref{cost_nu}), likelihoods (\ref{likelihood}),  
and topological energies (\ref{omega0})  are well defined.

\vskip 5mm

{\bf Acknowledgements.} 
This research has been partially supported by 
the FEDER /Ministerio de Ciencia, Innovaci\'on y Universidades - Agencia 
Estatal de Investigaci\'on grants No. MTM2017-84446-C2-1-R and 
PID2020-112796RB-C21. 
AC thanks G. Stadler  for nice discussions and useful suggestions.

\end{document}